\documentclass[12pt]{article}
\usepackage{amsmath}
\usepackage{amsfonts}
\usepackage{amssymb}
\usepackage{latexsym}
\usepackage{euscript}
\usepackage{a4wide}

\input xypic
\xyoption{all}

\def\frak{\mathfrak}

\def\SL{\mathrm{SL}}

\def\text{\textrm}

\def\GL{\mathrm{GL}}

\def\Spec{\mathop{\mathrm{Spec}}\nolimits}

\def\Qbar{\overline{\Q}}

\def\Gal{\mathrm{Gal}}
\def\F{\mathbf F}
\def\Fbar{\overline{\F}}

\def\nubar{\bar{\nu}}
\def\sigmabar{\bar{\sigma}}
\def\taubar{\bar{\tau}}
\def\pibar{\bar{\pi}}
\def\chibar{\overline{\chi}}
\def\rhobar{\overline{\rho}}
\def\rhouniv{\rho^{\mathrm{univ}}}
\def\rhomin{\rho^{\mathrm{min}}}
\def\qed{\hfill \square \ }
\def\Of{\mathcal O}
\def\E{{\mathcal E}}

\def\D{\mathcal D}

\def\Cal{\mathcal}
\def\Cl{\mathrm{Cl}}

\def\Q{\mathbf Q}
\def\Z{\mathbf Z}

\def\T{\mathbf T}
\def\twT{\widetilde{\T}}

\def\GL{\mathrm{GL}}

\def\Ext{\mathrm{Ext}}
\def\End{\mathrm{End}}

\def\E{\cal E}

\def\Kab{K^{ab}}
\def\et{\text{\'et}}

\def\co{0}
\def\Def{\mathrm{Def}}
\def\Defmin{\mathrm{Def}^{\mathrm{min}}}
\def\Defminprime{(\mathrm{Def}^{\mathrm{min}})'}
\def\Hom{\mathop{\mathrm{Hom}}\nolimits}
\def\kernel{\mathop{\mathrm{ker}}\nolimits}
\def\image{\mathop{\mathrm{im}}\nolimits}
\def\cokernel{\mathop{\mathrm{coker}}\nolimits}
\def\Frob{\mathop{\mathrm{Frob}}\nolimits}
\def\Trace{\mathop{\mathrm{Trace}}\nolimits}
\def\Det{\mathop{\mathrm{Det}}\nolimits}
\def\sdp{\medspace \times \kern -2.1pt 
 \vrule width0.4pt height 6pt depth -0.1pt \medspace}

\def\ilim{\displaystyle \lim_{\longrightarrow}}

\def\iso{\cong}

\newtheorem{theorem}{Theorem}[section]
\newtheorem{lemma}[theorem]{Lemma}

\newtheorem{cor}[theorem]{Corollary}

\newtheorem{prop}[theorem]{Proposition}

\def\Ker{\mathrm{Ker}}
\def\XYZ{XY\kern-.17em{Z}}
\def\fppf{f\kern-.02em{p}p\kern-.07em{f}}
\def\Kf{K_{\kern-.1em{f}}}
\def\ef{e_p}
\def\eftwo{e_2}
\def\efthree{e_3}
\def\gp{g_p}

\def\m{\mathfrak{p}}
\def\triv{1}
\def\Vbar{\overline{V}}
\def\Lbar{\overline{L}}
\def\Vmin{V^{\mathrm{min}}}
\def\Lmin{L^{\mathrm{min}}}
\def\Vuniv{V^{\mathrm{univ}}}
\def\Luniv{L^{\mathrm{univ}}}

\renewcommand\pmod{\mod}
\def\f1{f'}

\def\groups#1#2{\mathcal{C}(#1^n,#2)}

\begin{document}

\author{Frank Calegari\footnote{Supported in part by the American Institute of Mathematics.} \\ Matthew Emerton}
\title{On the Ramification of Hecke Algebras at Eisenstein Primes}
\maketitle

\section{Introduction}

Fix a prime $p$, and a modular residual representation
$\rhobar: G_{\Q} \rightarrow \GL_2(\Fbar_p)$.
Suppose $f$ is a normalised cuspidal Hecke eigenform of
some level $N$
and weight $k$ that gives rise to $\rhobar$, and let
$\Kf$ denote the extension of $\Q_p$ generated by the
$q$-expansion coefficients $a_n(f)$ of~$f$.
The field $\Kf$ is a finite extension
of $\Q_p$. What can one say about the extension $\Kf/\Q_p$?
  Buzzard \cite{buzz} has made the following conjecture:
if $N$ is fixed, and $k$ is allowed to vary,
 then the degree $[\Kf:\Q_p]$ is bounded
independently of $k$.

Little progress has been made on this conjecture so far;
indeed, very little seems to have been proven at all regarding
the degrees $[\Kf:\Q_p]$.  The  goal of this paper
is to consider a somewhat orthogonal question to Buzzard,
namely, to fix the weight and vary the level. Moreover, we only
 consider certain \emph{reducible} representations
$\rhobar$ that arise in Mazur's study of the Eisenstein
Ideal~\cite{eisenstein}. Our results
suggest  that the degrees $[\Kf:\Q_p]$ are, in fact,
arithmetically
significant. 

Suppose that $N \geq 5$ is prime, and that $p$ is a prime
that exactly divides the numerator of $(N-1)/12$.
Mazur \cite{eisenstein}
has shown that there is a weight two cuspform defined over $\Qbar_p$,
unique up to conjugation by $G_{\Q_p}$ (the Galois group of $\Qbar_p$
over $\Q_p$),
satisfying the congruence
\begin{equation}\label{congruence}
a_{\ell}(f) \equiv 1 + \ell \pmod{\m}\end{equation}
(where $\m$ is the maximal ideal in the ring of integers of $\Kf$,
and $\ell$ ranges over primes distinct from $N$).
It follows moreover from \cite{eisenstein} (Prop.~19.1, p.~140)
 that $\Kf$ is a \emph{totally ramified}
extension of $\Q_p$, and thus that the degree $[\Kf:\Q_p]$
is equal to the (absolute) ramification degree of $\Kf$.
Denote this ramification degree by $\ef$.

In this paper we prove the following theorem, in the case when
$p=2$.

\begin{theorem}\label{thm:main:p=2}
Suppose that $p = 2$ and that $N \equiv 9 \pmod 16$,
and let $f$ be 
a weight two eigenform on $\Gamma_0(N)$ satisfying the
congruence~(\ref{congruence}).
If $2^m$ is the largest power of $2$ dividing
the class number of the field $\Q(\sqrt{-N}),$
then $\eftwo = 2^{m-1}-1$.
\end{theorem}

When $p$ is odd, we establish the following less
definitive result.

\begin{theorem}\label{thm:main:p odd}
\label{theorem:oddp}
Suppose that $p$ is an odd prime exactly dividing
the numerator of $(N-1)/12$.
Let $f$ be 
a weight two eigenform on $\Gamma_0(N)$ satisfying the
congruence~(\ref{congruence}).
\begin{itemize}
\item[(i)] Suppose that $p=3$. $($Our hypothesis on $N$ thus becomes
$N \equiv 10 \text{ or } 19 \pmod{27})$. 
Then $\efthree = 1 $ if and only if the $3$-part of
the class group of $\Q(\sqrt{-3}, N^{1/3})$ is cyclic.
\item[(ii)] Suppose that $p \geq 5$. $($Our hypothesis on $N$
thus becomes $p \| N-1)$.  Then $\ef = 1 $ if 
the $p$-part of the class group of $\Q(N^{1/p})$ is cyclic.
\end{itemize}
\end{theorem}

The question of computing $\ef$ has been addressed previously,
in the paper \cite{merel} of Merel.  In this work, Merel
establishes a necessary and sufficient criterion for $\ef = 1$.
Merel's criterion for $\ef = 1$ is {\it not} expressed
in terms of class groups; rather, it is expressed in terms of
whether or not the
congruence class modulo $N$ of a certain explicit expression
is a $p$th power.

When $p = 2$, 
Merel, using classical results from algebraic number theory,
was able to reinterpret his explicit criterion for $\eftwo=1$
so as to prove that $\eftwo = 1$ if and only if $m = 2$.
(It is known that $m \geq 2$ if and only if $N \equiv 1 \pmod{8}.$)
Theorem \ref{thm:main:p=2} strengthens this result, by relating
the value of $\eftwo$ in all cases to the order of the 2-part of the class
group of $\Q(\sqrt{-N}).$

When $p$ is odd, Merel was not able to reinterpret his explicit
criterion in algebraic number theoretic terms.  However, combining
Merel's result with Theorem~\ref{thm:main:p=2}
(and the analogue of this theorem for more general $N$),
we obtain the following result. 

\begin{theorem}
\label{theorem:merel}
 Let $N \ge 5$ be prime.
\begin{itemize}
\item[(i)] Let $N \equiv 1 \pmod 9$. 
The $3$-part of the class group of $\Q(\sqrt{-3},N^{1/3})$
 is cyclic if and only if
$\left(\frac{N-1}{3} \right)!$ is not a cube modulo  $N$.
Equivalently, if we let $N = \pi \pibar$ denote the factorisation
of $N$ in $\Q(\sqrt{-3}),$ then the $3$-part of the class group
of $\Q(N^{1/3},\sqrt{-3})$ is cyclic if and only if the
$9$th power residue symbol $\left( \dfrac{\pi}{\pibar}\right)_9$
is non-trivial.
\footnote{The claimed equivalence follows from the
formula $\left(\left(\frac{N-1}{3}\right) !\right)^3 \equiv \pi
\pmod \pibar$, which was pointed out to us by Noam Elkies.
Ren\'e Schoof has told us that one can
prove part~(i) of Theorem~\ref{theorem:merel}
using class field theory. It is not apparent, however, that
(ii) can be proved in this way.}

Furthermore, if these equivalent conditions hold, then
the $3$-part of the class group of $\Q(N^{1/3})$ (which
a fortiori
is cyclic of order divisible by three) has order exactly three.

\item[(ii)] Let $p \ge 5$, and let $N \equiv 1 \pmod p$. 
If the $p$-part of the class group of $\Q(N^{1/p})$ is cyclic  then
$$\prod_{\ell=1}^{(N-1)/2} \ell^{\ell}$$
is not a $p$th power modulo $N$.
\end{itemize}
\end{theorem}

The proof of Theorems~\ref{thm:main:p=2} and \ref{thm:main:p odd}
depends on arguments using deformations of 
Galois representations.    Briefly, if $\T$ denotes
the completion of the Hecke algebra acting on weight two
modular forms on $\Gamma_0(N)$ at its $p$-Eisenstein
ideal, then we identify $\T$ with the universal deformation
ring for a certain deformation problem.   The theorems
are then proved by an explicit analysis of this deformation
problem over Artin $\F_p$-algebras.

It may be of independent interest to note that
our identification of $\T$ as a universal deformation
ring also allows us to recover
{\it all} the results of Mazur proved in the reference
\cite{eisenstein} regarding the
structure of $\T$ and the Eisenstein ideal:
for example, that $\T$ is monogenic over $\Z_p$
(and hence Gorenstein);
that the Eisenstein ideal is principal,
and is generated by $T_{\ell} - (1 + \ell)$
if and only if $\ell \neq N$ is a good prime;
and also that $T_N = 1$ in $\T$.

\

Let us now give a more detailed explanation of our
method.  For the moment, we 
relax our condition on $N$, assuming simply that
$N$ and $p$ are distinct primes.
We begin by defining a continuous representation
$\rhobar: G_{\Q} \rightarrow \GL_2(\F_p)$.
If $p$ is odd, we let
$$\rhobar = \left(\begin{matrix} \chibar_p & 0 \\ 0 & 1 \end{matrix}\right),$$
where $\chibar_p$ is the mod $p$ reduction of the cyclotomic
character.
If $p$ is even, we let
$$\rhobar = \left(\begin{matrix} 1 & \phi \\ 0 & 1 \end{matrix}\right),$$
where $\phi:G_{\Q}\rightarrow \F_2$ is the unique
$\F_2$-valued homomorphism
inducing an isomorphism $\Gal(\Q(\sqrt{-1})/\Q) \iso \F_2$.

Let $\Vbar$ denote the two dimensional vector space on which $\rhobar$
acts, and fix a line $\Lbar$ in $\Vbar$ that is {\it not} invariant
under $G_{\Q}$ (equivalently,
$G_{\Q_p}$).

If $A$ is an 
Artinian local ring with residue field $\F_p$,
consider the set of triples $(V,L,\rho)$, where $V$ is a 
free $A$-module, $L$ is a direct summand of $V$ that
is free of rank one over $A$, and $\rho$ is a continuous homomorphism
$G_{\Q} \rightarrow \GL(V)$, satisfying the 
following
 conditions:
\begin{enumerate}
\item The triple $(V,L,\rho)$ is a deformation of $(\Vbar, \Lbar,\rhobar)$.
\item The representation $\rho$ is unramified away from $p$ and $N$,
and is finite at $p$ (i.e.~$V$, regarded as a $G_{\Q_p}$-module, arises as
the generic fibre of a finite flat group scheme over~$\Z_p$).
\item The inertia subgroup at $N$ acts trivially on the submodule $L$ of $V$.
\item The determinant of $\rho$ is equal to the composition of
the cyclotomic character
$\chi_p: G_{\Q} \rightarrow \Z_p^{\times}$ with the natural map
$\Z_p^{\times} \rightarrow A^{\times}$.
\end{enumerate}

If we let $\Def(A)$ denote the collection of such triples
modulo strict equivalence, then $\Def$ defines a deformation
functor on the category of Artinian local rings $A$.

Note that the representation $\rhobar$ is reducible,
and is either the direct sum of two characters (if $p$ is odd)
or an extension of the trivial character by itself (if $p = 2$).
Nevertheless,
one has the following result.

\begin{prop}\label{representability}
The deformation functor $\Def$ is pro-representable by
a complete Noetherian local $\Z_p$-algebra $R$.
\end{prop}

The proposition follows directly from that fact that the only
endomorphisms of the triple $(\Vbar,\Lbar,\rhobar)$ are the scalars.
(The authors learned the idea of introducing a locally invariant line
to rigidify an otherwise unrepresentable
deformation problem from Mark Dickinson, who has applied it
to analyse the deformation theory of residually irreducible representations
that are ordinary, but not $p$-distinguished, locally at $p$.)

\

Having defined a universal deformation ring, we now introduce the
corresponding Hecke algebra.  As indicated above,
we let $\T$ denote the completion at its $p$-Eisenstein
ideal of the $\Z$-algebra of Hecke operators
acting on the space of all modular forms (i.e.~the cuspforms
together with the Eisenstein series) of level $\Gamma_0(N)$
and weight two.  (The $p$-Eisenstein ideal is
the maximal ideal in the Hecke algebra
generated by the elements $T_{\ell} - (1 + \ell)$
($\ell \neq N$), $T_N  - 1$, and $p$).

The following result relates $R$ and $\T$.

\begin{theorem}\label{R=T}
If $\rhouniv$ denotes the universal deformation of $\rhobar$
over the universal deformation ring $R$,  
then there is an isomorphism of $\Z_p$-algebras $R \iso \T$,
uniquely determined by the requirement that
the trace of Frobenius at $\ell$ under $\rhouniv$ $($for primes $\ell \neq p,N)$
maps to the Hecke operator $T_{\ell} \in \T$.
\end{theorem}

Let us now return to the setting of Theorems~\ref{thm:main:p=2}
and~\ref{thm:main:p odd}.
Thus we suppose again that $p$ exactly divides the numerator
of $(N-1)/12$,
and let $f$ be as in the statements of the theorems.
If $\Cal O$ denotes the ring of integers in $\Kf$,
and $\m$ its maximal ideal,
then the results of \cite{eisenstein} imply (taking into account
the congruence satisfied by $N$) that
the Hecke algebra $\T$ admits the following description:
$$\T
= \{(a,b) \in \Z_p \times \Cal O \, | \, a \bmod p = b \bmod \m \}.$$
From this description of $\T$, one easily computes that
$\T/p$ is isomorphic
to $\F_p[X]/X^{\ef+1}.$
Theorem~\ref{R=T} thus yields the following characterisation of $\ef$.

\begin{cor}\label{EPR}
The natural number $\ef$ is the largest integer $e$ for which
we may find a triple $(V,L,\rho)$ in $\Def(\F_p[X]/X^{e+1})$
such that the induced map 
$R \rightarrow \F_p[X]/X^{e+1}$ is surjective.
\end{cor}

Theorems~\ref{thm:main:p=2} and~\ref{thm:main:p odd} are
a consequence of this corollary,
together with an explicit analysis of the deformations
of $(\Vbar,\Lbar,\rhobar)$ over Artinian local rings of the
form $\F_p[X]/X^n$. 

If $p^2$ divides the numerator
of $(N-1)/12$, then the residually
Eisenstein cusp forms of level $N$ need not be
mutually conjugate.
However, one still has an isomorphism of the
form $\T/p = \F_p[x]/x^{\gp + 1}$, where $\gp+1$ denotes
the rank of $\T$ over $\Z_p$. (Thus $\gp$ is the
rank over $\Z_p$ of the cuspidal quotient of $\T$.)
In particular, the cuspidal Hecke algebra localized at the Eisenstein
prime is isomorphic to $\Z_p$ if and only if $\gp = 1$.
In this way
our analysis of deformations over $\F_p[X]/X^n$ suffices to
prove Theorem~\ref{theorem:merel}. More generally,
our paper can
be seen as providing a partial answer to Mazur's
question (\cite{eisenstein}, p.~140): ``\emph{Is there anything
general that can be said \ldots about $g_p$?}''.

The organisation of the paper is as follows.  In 
Section~\ref{sec:ext}
we develop some results about group schemes that will
be required in our study of the deformation functor $\Def$.
In Section~\ref{sec:RT} we prove Theorem~\ref{R=T},
using the numerical criterion of Wiles \cite{wiles}
(subsequently strengthened by Lenstra \cite{lenstra}). 
As in \cite{skiles}, we use the class field theory of cyclotomic
fields to obtain the required upper bound for the size of an appropriate
Galois cohomology group; the numerical criterion is then
established by comparing this upper bound with the congruence modulus
of the weight two Eisenstein series on $\Gamma_0(N)$ (which is known
by \cite{eisenstein} to equal the numerator of $(N-1)/12$). 
Finally in Sections~\ref{sec:explicit at 2} (respectively~\ref{sec:odd explicit})
we perform the analysis
necessary to deduce Theorem~\ref{thm:main:p=2} (respectively~\ref{thm:main:p odd})
from Corollary~\ref{EPR}.

Let us close this introduction by emphasising that the only
result of \cite{eisenstein} required for the proof
of Theorem~\ref{R=T}
is the computation of the congruence modulus
between the Eisenstein and cuspidal locus in
the Hecke algebra of weight two and level $N$.
(Namely, that this congruence modulus is equal
to the numerator of $(N-1)/12$.)
As remarked upon above, we are then able to deduce
all the results of \cite{eisenstein} regarding $\T$
and its quotient $\T^0$ from Theorem~\ref{R=T}.
The necessary arguments are presented
at the end of Section~\ref{sec:RT}.

%
%

\section{Some group scheme-theoretic calculations}
\label{sec:ext}
Let us fix a prime $p$, and a natural number $n$.
We begin by considering finite flat commutative groups
schemes of exponent $p^n$ that are extensions of
$\Z/p^n$ by $\mu_{p^n}$.
For any scheme $S$ we let $\groups{p}{S}$ denote
the category of
commutative finite flat group schemes of exponent $p^n$
over the base $S$,
and we write $\Ext^1_S(\text{--} , \text{--})$ to
denote the Yoneda Ext-bifunctor on the additive category 
$\groups{p}{S}$.

\begin{lemma}\label{inj-one}
The natural map $\Ext^1_{\Z_p}(\Z/p^n,\mu_{p^n})
\rightarrow \Ext^1_{\Q_p}(\Z/p^n,\mu_{p^n})$,
induced by restricting
to the generic fibre, is injective.
\end{lemma}

\begin{Proof}
Kummer theory identifies the map in the statement of the
lemma with the obviously injective map
$\Z_p^{\times}/(\Z_p^{\times})^{p^n} \rightarrow
\Q_p^{\times}/(\Q_p^{\times})^{p^n}. \qed$
\end{Proof}

\

%

If $p = 2$, we let $\Vmin_n$ denote the extension of
$\Z/2^n$ by $\mu_{2^n}$ in the category $\groups{2}{\Q}$
corresponding by Kummer theory to the element
$-1 \in \Q^{\times}/(\Q^{\times})^{2^n}$.
If $p$ is odd, we let $\Vmin_n$ denote the direct sum
$\Z/p^n \bigoplus \mu_{p^n}$ in the category $\groups{p}{\Q}$.
We may (and do) regard $\Vmin_n$ as an object of the category
of $G_{\Q}$-modules annihilated by $p^n$.

More explicitly, let $\chi_p$ denote the $p$-adic cyclotomic character.
Then if $p = 2,$ the $G_{\Q}$-module $\Vmin_n$ corresponds to the 
representation
$$\rhomin_n: G_{\Q} \rightarrow \GL_2(\Z/2^n)$$
given by
$$\sigma \mapsto \left( \begin{matrix} \chi_2(\sigma) &
(\chi_2(\sigma) - 1)/2 \\ 0 & 1 \end{matrix} \right) \pmod{2^n},$$
whilst if $p$ is odd, the $G_{\Q}$-module $\Vmin_n$ corresponds to
the representation
$$\rhomin_n: G_{\Q} \rightarrow \GL_2(\Z/p^n)$$
given by
$$\sigma \mapsto \left( \begin{matrix} \chi_p(\sigma)  & 0
\\ 0 & 1 \end{matrix}\right)\pmod p^n.$$
(Here we have denoted by $\sigma$ an element of $G_{\Q}$.)

\begin{prop}\label{unique prol}
For any natural number $M$, the $G_{\Q}$-module $\Vmin_n$
has a unique prolongation to an object of $\groups{p}{\Z[1/M]}$.
\end{prop}

\begin{Proof}
The Galois module $\Vmin_n$ is unramified away from $p$,
and so $\Vmin_n$ has a unique prolongation to an \'etale group
scheme over $\Z[1/Mp]$.
It thus suffices to show that $\Vmin_n$, regarded as
a $G_{\Q_p}$-module, has a unique prolongation to an object
of $\groups{p}{\Z_p}$.
If $p$ is odd, then this is a direct consequence
of \cite{fontaine}, Thm.~2.  Thus we assume
for the remainder of the proof that $p = 2$.
In this case, $\Vmin_n$ is defined to be the extension
of $\Z/2^n$ by $\mu_{2^n}$ corresponding to $-1 \in \Q_2^{\times}$.
Since $-1$ in fact lies in $\Z_2^{\times}$,
$\Vmin_n$ does prolong to a finite flat group scheme over $\Z_2$.
We must show that this prolongation is unique.

We begin with the case $n=1$.
Suppose that $G$ is a finite flat group
scheme over $\Z_2$ having $(\Vmin_1)_{/G_{\Q_2}}$ as its associated
Galois representation.  The scheme-theoretic closure of the fixed
line in $\Vmin_1$ yields an order two finite flat subgroup scheme $H$ of $G$.
Both $H$ and $G/H$ are thus finite flat group schemes of order two.
The results of \cite{OT} show that $\Z/2$ and $\mu_2$ are the only
group schemes of order 2 over $\Z_2$.
Thus $G$ is an extension of either $\Z/2$ or $\mu_2$ by
either $\Z/2$ or $\mu_2$.  Since neither $G$ nor its Cartier dual
are unramified (since $\Vmin_1$ is self-dual and ramified at 2),
we see that both $\Z/2$ and $\mu_2$ must appear.  Since $\Vmin_1$ is
a non-trivial $G_{\Q_2}$-module,
a consideration of the connected-\'etale exact sequence
attached to $G$ shows that in fact $G$ is an extension of $\Z/2$ by
$\mu_2$.  The fact that $G$ is determined uniquely by $\Vmin_1$
now follows from Lemma~\ref{inj-one}.

Now consider the case of $n$ arbitrary.  Let $G_n$ be an prolongation
of $\Vmin_n$ to a group scheme over $\Z_2$.  The Galois module
$\Vmin_n$ admits a filtration by submodules with successive
quotients isomorphic to $\Vmin_1$.
Taking scheme-theoretic closures, and appealing to the
conclusion of the preceding paragraph, we obtain a filtration
of $G_n$ by finite flat closed subgroup schemes, with successive quotients
isomorphic to the unique extension of $\Z/2$ by $\mu_2$ over $\Z_2$
with generic fibre isomorphic to $\Vmin_1$.
Consider the connected-\'etale sequence of $G_n$:
$$0 \rightarrow G^{\co} \rightarrow G \rightarrow G^{\et} \rightarrow 0.$$
We see that $G^{\co}$ is a successive extension of copies of 
$\mu_2$, and that $G^{\et}$ is a successive extension of copies
of $\Z/2$.  Such extensions (being respectively dual-to-\'etale,
or \'etale) are uniquely determined by their
corresponding Galois representations.  Thus $G^{\co}$ is in
fact isomorphic to $\mu_{2^n},$ whilst $G^{\et}$ is isomorphic to
$\Z/2^n$.  Thus $G_n$ is an extension of $\Z/2^n$ by $\mu_{2^n}$,
and Lemma~\ref{inj-one} shows that it is uniquely determined
by $\Vmin_n$.
$\qed$
\end{Proof}

\begin{lemma}\label{Exts}
Let $D_n$ denote the (uniquely determined, by Proposition~\ref{unique prol})
prolongation of $\Vmin_n$ to an object of $\groups{p}{\Z}$.
We have $\Ext^1_{\Z}(\Z/p^n,D_n) 
= 0$.
\end{lemma}

\begin{Proof}
Writing $D_n$ as an extension
of $\Z/p^n$ by $\mu_{p^n},$
we obtain the exact sequence of Yoneda Ext groups
$$\Ext^1_{\Z}(\Z/p^n,\mu_{p^n}) \rightarrow \Ext^1_{\Z}(\Z/p^n,D_n)
\rightarrow \Ext^1_{\Z}(\Z/p^n,\Z/p^n).$$
The third of these groups always vanishes, since $\Z$ has no non-trivial
\'etale covers.  If $p$ is odd, the first of these groups always vanishes.
If $p = 2$, then the first
of these groups has order two, with the non-trivial element corresponding
by Kummer theory to $-1 \in \Z^{\times}$.  Since $D_n$ is itself classified
by this same element when $p = 2$, we see that the first arrow vanishes in
all cases, and thus so does the middle group.
$\qed$
\end{Proof}

\begin{prop}\label{uniqueness}
Suppose that $M_{/\Q_p}$ is a $G_{\Q_p}$-module equipped with
a composition series whose subquotients are isomorphic to $\Vmin_1$,
that admits a prolongation
to a finite flat group scheme $M$ over $\Z_p$.
Then $M$ is unique up to unique isomorphism.
\end{prop}

\begin{Proof}
Let $M$ and $M'$ be two choices of a finite flat group scheme over $\Z_p$
prolonging $M_{/\Q_p}$.  The results of \cite{ray} show that we may find
a prolongation of $M_{/\Q_p}$ that maps (in the category of such prolongations)
to each of $M$ and $M'$.  Thus we
may assume we are given a map $M \rightarrow M'$ that induces the identity on
generic fibres.  By assumption we may find an embedding
$\Vmin_1 \subset M_{/\Q_p}$.
Passing to scheme theoretic closures in each of $M$ and $M'$,
and taking into account Lemma \ref{unique prol},
this prolongs to an embedding of $D_1$ into each of $M$ and $M'$,
so that the map $M \rightarrow M'$ restricts to the identity map
between these two copies of $D_1$.   
Replacing $M_{/\Q_p}$ by $M_{\Q_p}/\Vmin_1,$
$M$ by $M/D_1$, and $M'$ by $M'/D_1$, and arguing by induction on the order
of $M$, the proposition follows from the $5$-lemma (applied,
for example, in the category of sheaves on the $\fppf$ site
over $\Z_p$).
$\qed$
\end{Proof}

\begin{cor}\label{uniquenesscor}
Suppose that $A$ is an Artinian local ring with maximal ideal $\m$
and residue field $\F_p$, that $V$ is a free $A$-module of rank two,
and that $\rho:G_{\Q_p}
\rightarrow \GL(V)$ is a deformation of $(\Vmin_1)_{/ \Q_p}$ that is 
finite flat at $p$.
Then there is a unique up to unique isomorphism finite flat group scheme
$M$ over $\Z_p$
whose generic fibre equals $V$.  Furthermore, the $A$-action on
$V$ prolongs to an $A$-action on $M$, and the connected-\'etale sequence
of $M$ gives rise to a two-step filtration of $V$ by free $A$-submodules
of rank one.
\end{cor}

\begin{Proof}
If we choose a Jordan-H\"older filtration of $A$ as a module over itself,
then this induces a composition series on $V$ with successive quotients
isomorphic to $\Vmin_1$.  Thus we are in the situation of the
preceding Proposition, and the uniqueness of $M$ follows.  
In particular $M$ is the maximal prolongation of $V,$ and by
functoriality, the $A$-action on $V$ prolongs to an $A$-action
on $M$. 
on $M$.  (More precisely, just as in the discussion of
\cite{ray}, p.~265, we see that the automorphisms of $V$ induced
by the group of units $A^{\times}$ extend to automorphisms of $M$.
Since $A$ is generated as a ring by $A^{\times}$, we conclude that in
fact the $A$-action on $V$ extends to an $A$-action on $M$.)
Finally, let $0 \rightarrow M^{\co} \rightarrow M
\rightarrow M^{\et} \rightarrow 0$ be the connected-\'etale sequence of $M$.
The functorial nature of its construction implies that it is an exact sequence
of $A$-submodules of $M$.
Thus the exact sequence $0 \rightarrow M^{\co}_{/\Q_p} \rightarrow V
\rightarrow M^{\et}_{/\Q_p} \rightarrow 0$ obtained by restricting to
$\Q_p$ yields a two-step filtration of $V$ by $A$-submodules.
The formation of this filtration is clearly functorial in $A$. 
Thus if we tensor
$M$ with $A/\m$ over $A$, and take into account
that $A/\m \otimes_A M = D_1,$ we find that
$A/\m \otimes_A M^{\co} = D_1^{\co} = \mu_p$
and that
$A/\m \otimes_A M^{\et} = D_1^{\et} = \Z/p.$
Thus each of $M^{\co}_{/\Q_p}$ and $M^{\et}_{/\Q_p}$ are cyclic
$A$-modules.  Since $M_{/\Q_p}$ is free of rank two over $A$,
they must both be free $A$-modules of rank one, as claimed.
$\qed$
\end{Proof}

\begin{prop}\label{ordinary}
Let $A$ be an Artinian local $\Z/p^n$-algebra
with maximal ideal $\m$ and residue field $\F_p$. 
If $V$ is a free $A$-module of rank two and $\rho:G_{\Q_p}
\rightarrow \GL(V)$ is a deformation of $(\Vmin_1)_{/\Q_p}$
that is finite at $p$,
then the coinvariants of $V$ with respect to the
inertia group $I_p$ are free of rank one over $A$.
\end{prop}

\begin{Proof}
The preceding corollary shows that $V$ admits a two-step
filtration, with rank one free quotients, corresponding
to the connected-\'etale sequence of the prolongation of
$V$ to a group scheme over $\Z_p$.  In particular,
the inertial coinvariants $V_{I_p}$ admit a surjection
onto a free $A$-module of rank one.
On the other hand, if $\m$ is the maximal ideal
of $A$, then
$(V_{I_p})/\m = (V/\m)_{I_p} = (\Vmin_1)_{I_p}.$
This latter space is directly checked to be one dimensional
over $\F_p$, implying that $V_{I_p}$ is a cyclic $A$-module.  
Altogether, we find that $V_{I_p}$ is free of rank one over
$A$, as claimed.
$\qed$
\end{Proof}

\begin{prop}\label{unram def}
Let $A$ be an Artinian local $\Z/p^n$-algebra
with maximal ideal $\m$ and residue field $\F_p$. 
If $V$ is a free $A$-module of rank two and $\rho:G_{\Q}
\rightarrow \GL(V)$ is a deformation of $\Vmin_1$ that is unramified away
from $p$ and finite at $p$, 
then there is an isomorphism
$V \cong A\otimes_{\Z/p^n} \Vmin_n.$
\end{prop}

\begin{Proof}
Corollary~\ref{uniquenesscor} shows that $V$ prolongs to a finite flat
$A$-module scheme $M$ over $\Spec \Z$.  If we choose a Jordan-H\"older
filtration of $A$ as an $A$-module, then this gives rise to a corresponding
filtration of $V$, with successive quotients isomorphic to $\Vmin_1$.
Passing to the scheme-theoretic closure in $M$, and taking into
account Proposition~\ref{unique prol}, we obtain
a filtration of $M$ by closed subgroup schemes, with successive quotients
isomorphic to $D_1$.  

If $p$ is odd, then Proposition~I.4.5 of~\cite{eisenstein}
shows that $M$ is (in a canonical way) the product of a constant group scheme
and a $\mu$-type group (that is, the Cartier dual of a constant group scheme).
Each of these subgroups is then an $A$-module scheme, and we easily conclude
that $V \cong A\otimes_{\Z/p^n} \Vmin_n.$

If $p=2$, then Propositions~I.2.1 and~I.3.1 of \cite{eisenstein}
show that $M$ is the extension of a constant group scheme
by a $\mu$-type group.  Again, each of these groups is seen to be
an $A$-module scheme, and we easily conclude that $M$ is in fact an extension
of the constant $A$-module scheme $A$ by the $\mu$-type $A$-module scheme
$A\otimes_{\Z/2^n}\mu_{2^n}$.
The group of all such extensions is classified
by
$$H^1(\Spec \Z, A\otimes_{\Z/2^n} \mu_{2^n})\iso A\otimes_{\Z/2^n}
\Z^{\times}/(\Z^{\times})^{2^n} \iso A/\frak{p} \otimes_{\F_2} \{\pm 1\}.$$
Since $V$ is a deformation of $\Vmin_1$, which corresponds by Kummer
theory to
the non-trivial element of $\{\pm 1\}$, we see that in
fact $M$ is classified by the non-trivial element of
$A/\frak{p} \otimes_{\F_2}\{\pm1\},$
and thus that $M \iso A\otimes_{\Z/2^n} D_n,$
and hence that $V\iso A\otimes_{\Z/2^n} \Vmin_n.$
$\qed$
\end{Proof}

\

We leave it to the reader to verify
the following lemma.

\begin{lemma}\label{minendos}
If $A$ is an Artinian local $\Z/p^n$-algebra
with maximal ideal $\m$,
then the ring of Galois equivariant endomorphisms of
$A\otimes_{\Z/p^n}\Vmin$ admits
the following description:
\begin{itemize}
\item[(i)] If $p = 2,$ then $\End_{A[G_{\Q}]}(A\otimes_{\Z/2^n}\Vmin_n)
=\{ \left(
\begin{matrix} a & b \\ 0 & a - 2b \end{matrix} \right) \, | \, a,b \in A\}.$
\item[(ii)] If $p$ is odd, then $\End_{A[G_{\Q}]}(A\otimes_{\Z/p^n}\Vmin_n)
=\{ \left(
\begin{matrix} a & 0 \\ 0 & d \end{matrix} \right) \, | \, a,d \in A\}.$
\end{itemize}
\end{lemma}

%

\section{Proving that $R=\T$}
\label{sec:RT}

In this section we prove Theorem \ref{R=T}.  We let $\Def$ denote
the deformation problem described in the introduction. 
We employ the technique introduced in~\cite{wiles}: namely,
we first consider a minimal deformation $\rhomin$ of $(\Vbar,\Lbar,\rhobar)$
over $\Z_p$, and then verify the numerical criterion of~\cite{wiles}.

Let us define the minimal deformation problem $\Defmin$, as
the subfunctor of $\Def$ consisting of those deformations 
of $(\Vbar,\Lbar,\rhobar)$ that are unramified away from $p$.
Let us also define $\Defminprime$ to be the functor that
classifies all deformations of $(\Vbar,\rhobar)$ that are unramified 
away from $p$ and finite at $p$.  Forgetting the $I_N$-fixed line
$L$ gives a natural transformation $\Defmin \rightarrow \Defminprime.$

Let us now define the triple $(\Vmin,\Lmin,\rhomin)$.
We take $\Vmin = \Z_p\bigoplus \Z_p$.
If $p = 2$, then we 
let $\rhomin$ denote the representation
$$\sigma \mapsto \left( \begin{matrix} \chi_2(\sigma) &
(\chi_2(\sigma) - 1)/2 \\ 0 & 1 \end{matrix} \right)$$
$($here $\sigma$ denotes an element of $G_{\Q})$,
while if $p$ is odd, we let $\rhomin$ denote the direct sum
of $\chi_p$ (the $p$-adic cyclotomic character) and $\triv$ (the trivial
character).
In each case,
the pair $(\Vmin,\rhomin)$ is certainly a lifting of $(\Vbar,\rhobar)$.
We take $\Lmin$ to be any free of rank one $\Z_p$-submodule of $\Vmin$
lifting the line $\Lbar$ in $\Vbar$.  

Note that 
for any natural number $n$, we have $\Vmin/p^n = \Vmin_n$ (the Galois module
introduced in the preceding section).

\begin{prop}\label{mindefiso} 
 The natural transformation $\Defmin \rightarrow \Defminprime$
is an isomorphism of functors. Moreover,
the deformation functor $\Defmin$ is pro-represented by 
$(\Vmin,\Lmin,\rhomin)$ in $\Defmin(\Z_p).$
\end{prop}

\begin{Proof}
Let $A$ be an Artinian local $\Z_p$-algebra, and
let $(V,\rho)$ be an object of $\Defminprime(A)$.
Proposition~\ref{unram def} shows that there is
an isomorphism $V \iso A\otimes_{\Z_p} \Vmin.$
The explicit description of the endomorphisms of
$A\otimes_{\Z_p}\Vmin$ provided by Lemma~\ref{minendos}
shows that we may furthermore choose this isomorphism
so that it is strict.  Thus we see that
$\Defminprime$ is pro-represented by $\Z_p$, with
$(\Vmin,\rhomin)$ as universal object.

Now suppose that $(V,L,\rho)$ is an object
of $\Defmin(A).$  Using Lemma~\ref{minendos}
again, we see that we may choose the strict endomorphism
$V \iso A\otimes_{\Z_p} \Vmin$ of the preceding paragraph
in such a way that
$L$ is identified with $A\otimes_{\Z_p} \Lmin.$
Thus $\Z_p$ also pro-represents $\Defmin,$ with universal
object $(\Vmin,\Lmin,\rhomin).$  This establishes the proposition.
$\qed$
\end{Proof}

\

Note that the preceding lemma implies in particular that
the class of $(\Vmin,\Lmin,\rhomin)$ in $\Def(\Z_p)$ is independent
of the choice of $\Lmin$ (provided that it lifts $\Lbar$).
Of course, this is easily checked directly, using the description
of $\End_{G_{\Q}}(\Vmin)$ afforded by Lemma~\ref{minendos}.

\

Let $R$ denote the universal deformation ring that pro-represents
the functor $\Def$, and
let $(\Vuniv,\Luniv,\rhouniv)$ denote the
universal deformation of $(\Vbar,\Lbar,\rhobar)$ over~$R$.
Corresponding to $(\Vmin,\Lmin,\rhomin)$ there is
a homomorphism $R \rightarrow \Z_p$ of $\Z_p$-algebras.
We let $I$ denote the kernel of this homomorphism.
The following more explicit description of $I$ will be useful.

\begin{prop}\label{I described}  
If $S$ is any finite set of primes containing $p$ and $N$,
then $I$ is generated by the set
$$\{1 + \ell - \Trace(\rhouniv(\Frob_{\ell})) \, | \, \ell \not\in S\}.$$
\end{prop}

\begin{Proof}
Let $I_S$ denote the ideal generated by the stated set.
Clearly $I_S \subset I.$  We will show that the Galois representation
$G_{\Q} \rightarrow \GL_2(R/I_S)$ obtained by reducing $\rhouniv$
modulo $I_S$ is unramified at $N$.  It will follow from 
Proposition~\ref{mindefiso} that $I \subset I_S$, and the proposition
will be proved.  The argument is a variation of that used to prove
Prop.~2.1 of~\cite{skiles}.

Suppose first that $p$ is odd.  Let us choose a basis for $\Vuniv$,
and write 
$$\rhouniv(\sigma) = \left( \begin{matrix} a({\sigma}) & b(\sigma) \\
c(\sigma) & d(\sigma)\end{matrix} \right),$$
for $\sigma \in G_{\Q}$.  We may assume that if $c \in G_{\Q}$ denotes
complex conjugation, then
$$\rhouniv(c) = \left( \begin{matrix} -1 & 0 \\ 0 & 1 \end{matrix}
\right).$$
We find that
$$a(\sigma) = \dfrac{1}{2}\left( \Trace(\rhouniv(\sigma))
- \Trace(\rhouniv(c \sigma)) \right),$$
and that
$$d(\sigma) = \dfrac{1}{2}\left( \Trace(\rhouniv( \sigma))
+ \Trace(\rhouniv(c \sigma))\right).$$
Thus 
$$a(\sigma) \equiv \chi_p(\sigma) \pmod{I_S},$$
whilst 
$$d(\sigma) \equiv 1 \pmod{I_S}.$$
In particular, if $\sigma \in I_N,$ then
$$\rhouniv(\sigma) \equiv \left( \begin{matrix} 1 & b(\sigma) \\
c(\sigma) & 1 \end{matrix} \right) \pmod{I_S}.$$
The universal $I_N$-fixed line is spanned by a vector
of the form $(1,x),$ where $x \in R^{\times}$.
We conclude that if $\sigma \in I_N$ then
$$(1 + b(\sigma) x, c(\sigma) + x) \equiv (1,x) \pmod{I_S},$$
and thus that
$$b(\sigma) \equiv c(\sigma) \equiv 0 \pmod{I_S}.$$
This implies that $\rhouniv \pmod{I_S}$ is unramified at $N$,
as required.

Consider now the case $p = 2$.
Again, we write
$$\rhouniv(\sigma) = \left( \begin{matrix} a({\sigma}) & b(\sigma) \\
c(\sigma) & d(\sigma)\end{matrix} \right),$$
for $\sigma \in G_{\Q}$.  We may assume that if $c \in G_{\Q}$ denotes
complex conjugation, then
$$\rhouniv(c) = \left( \begin{matrix} -1 & -1 \\ 0 & 1 \end{matrix}\right).$$
We may also assume that the universal $I_N$-fixed line is spanned by
the vector $(0,1)$.
By considering $\Trace(\rhouniv(c \sigma)),$ for $\sigma \in G_{\Q},$ 
we find that
$$-a(\sigma) - c(\sigma) + d(\sigma) \equiv 1 - \chi_2(\sigma) \pmod{I_S}.$$
If $\sigma \in I_N,$ then since $\sigma$ fixes $(0,1),$ we find
that
$$b(\sigma) \equiv 0, \quad d(\sigma) \equiv 1 \pmod{I_S}.$$
The preceding equations, the fact that $\det \rho^{univ} = \chi_2,$
and the fact that $\chi_2(\sigma) = 1$
for $\sigma \in I_N,$ imply that also
$$a(\sigma) \equiv 1, \quad c(\sigma) \equiv 0 \pmod{I_S}.$$
Altogether, we conclude that $\rhouniv \pmod{I_S}$ is unramified
at $N$, as required.
$\qed$
\end{Proof}

\

The preceding result has the following important corollary.

\begin{cor}\label{traces}
If $S$ is any finite set of primes containing $p$ and $N$,
then the complete local $\Z_p$-algebra $R$ is topologically
generated by the elements $\Trace(\rhouniv(\Frob_{\ell}))$,
for $\ell \not\in S$.
\end{cor}

\begin{Proof}
This follows immediately from the description of
$I$ provided by Proposition~\ref{I described}, the fact that $R$
is $I$-adically complete, and the fact that $R/I \iso \Z_p$.
$\qed$
\end{Proof}

\

We now compute the order of $I/I^2$, which is one of the
two ingredients we will eventually use in our verification
of the Wiles-Lenstra numerical criterion.

\begin{theorem}\label{cohomology calculation} 
The order of $I/I^2$ $($which is a power of $p)$ divides $(N^2-1)/24$.
\end{theorem}

\begin{Proof}
Let $n$ be a natural number, and let $(\Vmin_n,\Lmin_n,\rhomin_n)$
denote the reduction modulo $p^n$ of $(\Vmin,\Lmin,\rhomin).$
We consider extensions of Galois modules
$$0 \rightarrow (\Vmin_n,\Lmin_n)
\rightarrow (E,F) \rightarrow (\Vmin_n,\Lmin_n) \rightarrow 0;$$
here the 
notion indicates that $E$ is a $G_{\Q}$-module that extends $\Vmin_n$
by itself, and that $F$ is a submodule of $E$ (not assumed to
be Galois invariant) providing an extension
of $\Lmin_n$ by itself.
We let $A_n$ denote the group of isomorphism classes of such extensions
for which
$E$ is annihilated by $p^n$, is unramified away from $p$ and $N$,
is finite at $p$, and is Cartier self-dual as an extension of
$\Vmin$ by itself, and for which $F$ is fixed (element-wise)
by the inertia group $I_N$.
%
%
%
The usual identification of the relative tangent space
to a deformation functor with an appropriate Ext-group in an
appropriate category of Galois modules shows that
$\Hom(I/I^2, \Q_p/\Z_p) \iso \ilim A_n.$

\begin{lemma}\label{one} If $(E,F)$ is an object of $A_n$ for which $E$ is
a trivial extension, then the pair $(E,F)$ is also a trivial extension.
\end{lemma}

\begin{Proof}
Let us remind the reader that
if $E$ is the trivial extension of $\Vmin_n$ by itself,
then the automorphisms of $E$ are of the form
$\left( \begin{matrix} \mathrm{Id} & A \\ 0 & \mathrm{Id} \end{matrix}
\right),$
where $A$ is an element of $\End_{G_{\Q}}(\Vmin).$
This being said, the lemma is easily checked using 
Lemma~\ref{minendos}.

Alternatively, we may appeal to  Proposition~\ref{mindefiso}.
Since $E$ is assumed to be a trivial extension, it is in
particular unramified at $N$, and thus corresponds to a
deformation for the subproblem $\Defmin$ of $\Def$.  The triviality of
$E$ implies that this deformation is trivial, when regarded as an
deformation for the problem $\Defminprime.$  Since $\Defmin$ maps isomorphically
to $\Defminprime$, we obtain the assertion of the lemma.
$\qed$
\end{Proof}

\

If $(E,F)$ is an object of $A_n,$ then since $\Z/p^n$ (with
the trivial $G_{\Q}$-action) is a quotient of $\Vmin_n$, whilst
$\mu_{p^n}$ (with its natural $G_{\Q}$-action) is a submodule
of $\Vmin_n$, the extension $E$ determines an extension
$E'$ of $G_{\Q}$ modules
\begin{equation}\label{extension}
 0 \rightarrow \Z/p^n \rightarrow E' \rightarrow \mu_{p^n}
\rightarrow 0.\end{equation}

\begin{lemma}\label{two} If $(E,F)$ is an object of $A_n$ for
which the extension $E'$ is trivial, then $E$ is also a trivial
extension.
\end{lemma}

\begin{Proof}
We let $D_n$ denote the (unique, by Proposition~\ref{unique prol})
prolongation of $\Vmin_n$ to a finite flat group scheme over $\Z[1/N]$.
Proposition~\ref{uniqueness} shows that $E$ has a unique prolongation
to a finite flat group scheme $\E$ over $\Z[1/N]$, that provides
an extension of $D_n$ by itself. 
We let $D_n^{(1)}$ denote the copy of $D_n$ that appears
as a submodule of $\E$, and let $D_n^{(2)}$ denote the copy
of $D_n$ that appears as a quotient.  Also, we let $\mu_{p^n}^{(i)}$
(respectively $(\Z/p^n)^{(i)}$) denote the copy of $\mu_{p^n}$
(respectively $\Z/p^n$) that appears as a subgroup scheme (respectively
a quotient group scheme) of $D_n^{(i)}$, for $i = 1,2$.

The quotient ${\E}/\mu_{p^n}^{(1)}$ is an extension of $D_n^{(2)}$ by
$(\Z/p^n)^{(1)}$.  Thus it yields a class
$e \in \Ext_{\Z[1/N]}(D_n^{(2)},(\Z/p^n)^{(1)})$. 
This latter group sits in the exact sequence
$$\Ext_{\Z[1/N]}((\Z/p^n)^{(2)},(\Z/p^n)^{(1)}) \rightarrow
\Ext_{\Z[1/N]}(D_n^{(2)},(\Z/p^n)^{(1)}) \rightarrow
\Ext_{\Z[1/N]}(\mu_{p^n}^{(2)},(\Z/p^n)^{(1)}).$$
By assumption the image of $e$ under the second arrow vanishes,
and thus there is a class
$e' \in \Ext_{\Z[1/N]}((\Z/p^n)^{(2)},(\Z/p^n)^{(1)})$
that maps to $e$ under the first arrow.
We can construct such an extension class $e'$ concretely as follows:
we may choose a lift of $\mu_{p^n}^{(2)}$ to a subgroup scheme $\mu$
of ${\E}/\mu_{p^n}^{(1)}$.  The quotient
$({\E}/\mu_{p^n}^{(1)})/\mu$
is then an extension of $(\Z/p^n)^{(2)}$ by $(\Z/p^n)^{(1)}$,
which gives a realisation of a class $e'$ mapping to $e$.
Our assumption on the submodule $F$ of $E$ implies that
it maps surjectively onto
$(E/\mu_{p^n}^{(1)})/\mu$
(the generic fibre of
$({\E}/\mu_{p^n}^{(1)})/\mu$),
and thus that the action of inertia at $N$ on
$(E/\mu_{p^n}^{(1)})/\mu$ is trivial.  Thus 
$({\E}/\mu_{p^n}^{(1)})/\mu$ has a prolongation
to a finite flat group scheme over $\Z$, yielding an extension of
$\Z/p^n$ by itself.  There are no such non-trivial extensions
that are finite flat over $\Z$, and thus the extension class $e'$ is
trivial.  Hence the extension class $e$ is also trivial, and so
${\E}/\mu_{p^n}^{(1)}$ is a split extension of $D_n^{(2)}$ by
$(\Z/p^n)^{(1)}$. 

If ${\E}'$ denotes the preimage in ${\E}$ of the subgroup $\mu_{p^n}^{(2)}
\subset D_n^{(2)}$, then (since ${\E}$ is Cartier self-dual),
we find that ${\E}'$ is Cartier dual to ${\E}/\mu_{p^n}^{(1)}$.
The result of the preceding paragraph thus shows that
${\E}'$ is a split extension of $\mu_{p^n}^{(2)}$ by $D_n^{(1)}$.
Consider the exact sequence
$$\Ext_{\Z[1/N]}((\Z/p^n)^{(2)},D_n^{(1)}) \rightarrow
\Ext_{\Z[1/N]}(D_n^{(2)},D_n^{(1)})\rightarrow
\Ext_{\Z[1/N]}(\mu_{p^n}^{(2)},D_n^{(1)}).$$
If $e''$ denotes the class of ${\E}$ in the middle group,
then we have just seen that its image under the second arrow
vanishes.  Thus we may find a class
$e''' \in \Ext_{\Z[1/N]}((\Z/p^n)^{(2)},D_n^{(1)})$
mapping to $e''$ under the first arrow.  
We can construct such a class $e'''$ concretely as follows:
lift $\mu_{p^n}^{(2)}$ to a subgroup scheme $\mu'$ of ${\E}'$.
The quotient ${\E}/\mu'$ then provides an extension
of $(\Z/p^n)^{(2)}$ by $D_n^{(1)}$ whose extension class
$e'''$ maps to $e''$.  Our assumption on $F$ implies
that its image in $E/\mu'$ (the generic fibre of ${\E}/\mu'$)
has non-zero image in $(\Z/p^n)^{(2)}$, and thus that
inertia at $N$ acts trivially on $E/\mu'$,
and so ${\E}/\mu'$ has a prolongation to a finite flat group
scheme over $\Z$ that extends $\Z/p^n$ by $D_n$.
Lemma~\ref{Exts} shows that any such extension is split,
and thus that $e'''$ vanishes.  Consequently $e''$ also
vanishes, and so $E$ is a split extension, as claimed.
$\qed$
\end{Proof}

\

Let $B_n$ denote the set of extensions
of $G_{\Q}$-modules of the form~(\ref{extension}) that
are unramified away from $p$ and $N$, 
and that prolong over $\Z_p$ to an extension of the finite flat group scheme
$\mu_{p^n}$ by the finite flat group scheme $\Z/p^n$.

\begin{lemma}\label{three} If $n$ is sufficiently large,
then the natural map $B_n \rightarrow B_{n+1}$ is an isomorphism,
and each side has order at most the $p$-power part of $(N^2-1)/24$.
\end{lemma}

\begin{Proof}
Let $\Sigma = \{p,N,\infty\},$ and let $G_{\Sigma}$ denote the
Galois group of the maximal extension of $\Q$ in $\Qbar$
unramified away from the elements of $\Sigma$.
Extensions of the form~(\ref{extension}) are classified by
the Galois cohomology group
$H^1(G_{\Sigma}, \mu_{p^n}^{\otimes -1}).$
If such an extension prolongs to an extension of finite flat
groups over $\Z_p$, then it is in fact trivial locally at $p$,
since the connected group scheme $\mu_{p^n}$ cannot have
a non-trivial extension over the \'etale group scheme $\Z/p^n$.
Thus $B_n$ 
is equal to the kernel of the natural map
$$H^1(G_{\Sigma}, \mu_{p^n}^{\otimes -1})\rightarrow
H^1(G_{\Q_p}, \mu_{p^n}^{\otimes -1}).$$

Let $K_n$ denote the extension of $\Q$ obtained by adjoining 
all $p^n$th roots of unity in $\Qbar$.  Let $H$ denote the
normal subgroup of $G_{\Sigma}$ which fixes $K_n$; the quotient
$G_{\Sigma}/H$ is naturally isomorphic to $(\Z/p^n)^{\times}$.
The prime $p$ is totally ramified in $K_n$. Thus, if $\pi$ denotes
the unique prime of $K_n$ lying over $p$, the quotient
$G_{\Q_p}/G_{K_{n,\pi}}$ also maps isomorphically to $(\Z/p^n)^{\times}$.
The inflation-restriction exact sequence gives a diagram
$$\xymatrix{
0 \ar[r] & H^1((\Z/p^n)^{\times},\mu_{p^n}^{\otimes {-1}}) \ar[r]\ar@{=}[d] &
H^1(G_{\Sigma},\mu_{p^n}^{\otimes{-1}}) \ar[r]\ar[d] &
H^1(H,\mu_{p^n}^{\otimes{-1}})^{(\Z/p^n)^{\times}} \ar[d] \\
0 \ar[r] & H^1((\Z/p^n)^{\times},\mu_{p^n}^{\otimes {-1}})\ar[r] &
H^1(G_{\Q_p},\mu_{p^n}^{\otimes{-1}}) \ar[r] &
H^1(G_{K_{n,\pi}},\mu_{p^n}^{\otimes{-1}})^{(\Z/p^n)^{\times}}.}$$
Taking into account the discussion of the preceding paragraph,
this diagram in turn induces an injection
$$B_n \hookrightarrow 
\kernel\,(
H^1(H,\mu_{p^n}^{\otimes{-1}})^{(\Z/p^n)^{\times}} \rightarrow \\
H^1(G_{K_{n,\pi}},\mu_{p^n}^{\otimes{-1}})^{(\Z/p^n)^{\times}}).$$
Since $H$ acts trivially on $\mu_{p^n}^{\otimes -1}$,
there is an isomorphism
$$H^1(H,\mu_{p^n}^{\otimes{-1}})^{(\Z/p^n)^{\times}} \iso
\Hom_{(\Z/p^n)^{\times}}(H,\mu_{p^n}^{\otimes -1}).$$
Thus $B_n$ injects into the subgroup of
$\Hom_{(\Z/p^n)^{\times}}(H,\mu_{p^n}^{\otimes -1})$
consisting of homomorphisms that are trivial on
$G_{K_{n,\pi}}.$

Any element of $\Hom_{(\Z/p^n)^{\times}}(H,\mu_{p^n}^{\otimes -1})$
that is trivial on $G_{K_{n,\pi}}$
factors through the Galois group $\Gal(L_n/K_n)$, where $L_n$
is the extension of $K_n$ defined in the statement of the following
lemma.
Lemma~\ref{three} is now seen to follow from the conclusion of that lemma.
$\qed$
\end{Proof}

\begin{lemma}\label{four} Let $L_n$ denote the maximal abelian extension 
of $K_n$ of degree dividing $p^n$ that is unramified away from $N$,
in which the prime lying over $p$ splits completely,
and on whose Galois group $(\Z/p^n)^{\times}
= \Gal(K_n/\Q)$ acts via $\chi_p^{-1}$.  Then $L_n$ is a cyclic
extension of $K_N$, and the degree
$[L_n:K_n]$ divides the $p$-power part of $(N^2 - 1)/24$.
\end{lemma}

\begin{Proof}
Let $\zeta$
be a choice of primitive $p^n$th root of unity.
If $\Of_n$ denotes
the ring of integers in $K_n$, then $1-\zeta$ generates the
unique prime ideal of $\Of_n$ lying above $p$.
Let $((\Of_n/N)^{\times}/p^n)_{(-1)}$ denote the maximal quotient of 
$(\Of_n/N)^{\times}/p^n$ on which $\Gal(K_n/\Q)$ acts via $\chi_p^{-1}$.
Since Herbrand's criterion shows that the $\chi_p^{-1}$-eigenspace
in the $p$-part of the class group of $K_n$ vanishes,
global class field theory shows that the Galois
group of $L_n/K_n$ is equal to
the cokernel of the composite
$$\Of_n^{\times}[1-\zeta,(1-\zeta)^{-1}] \rightarrow (\Of_n/N)^{\times}
\rightarrow ((\Of_n/N)^{\times}/p^n)_{(-1)}.$$

Fix a prime $\frak n$ of $K_n$ lying over $N$.   The inclusion
$(\Of_n/\frak n)^{\times} \hookrightarrow (\Of_n/N)^{\times}$
induces an isomorphism
\begin{equation}\label{induction formula}
(\Of_n/\frak n)^{\times}/(p^n,N^2-1)
\iso ((\Of_n/N)^{\times}/p^n)_{(-1)}.
\end{equation}
Since $(\Of_n/\frak n)^{\times}$ is a cyclic group,
we conclude that $((\Of_n/N)^{\times}/p^n)_{(-1)}$
has order bounded by the $p$-part of $N^2-1$.
Thus if $p \geq 5$ the lemma is proved.

We now perform a more refined analysis, which will prove
the lemma in the remaining cases (i.e.~$p = 2$ or 3).
A simple computation shows that under the isomorphism~(\ref{induction formula}),
the subgroup of 
$$((\Of_n/N)^{\times}/p^n)_{(-1)}$$ generated by $(1-\zeta)$
corresponds to the subgroup of $(\Of_n/\frak n)^{\times}/(p^n,N^2-1)$
generated by
$$\prod_{a \in (\Z/p^n)^{\times}/\langle N \rangle } (1 - \zeta^a)^a$$
(where $\langle N \rangle$ denotes the cyclic subgroup of $(\Z/p^n)^{\times}$
generated by $N \mod p^n$).

Suppose first that $p$ is odd, so that the order of $(\Z/p)/\langle N \rangle$
is prime to $p$.  Let $c$ denote this order.  Also, write
$N = \omega(N) N_1$ in $\Z_p$, where $\omega(N)$ is the Teichm\"uller lift
and $N_1$ is a 1-unit.  If $p^f$ denotes the exact power of $p$ dividing
$N^2-1$, and $p^{\f1}$ denotes
the exact power of $p$ dividing $N_1 - 1$,
then $\f1 \geq f,$ with equality if $p = 3$.  Let us assume that
$n \geq \f1$, so that
$(\Of_n/\frak n)^{\times}/(p^n,N^2-1)$ is cyclic of order $p^f,$
generated by the image of $\zeta$, or of $-\zeta$.

Since $2c$ is prime to $p,$
the subgroup of $(\Of_n/\frak n)^{\times}/(p^n,N^2-1)$
generated by
$$\prod_{a \in (\Z/p^n)^{\times}/\langle N \rangle } (1 - \zeta^a)^a$$
coincides with the subgroup generated by
$$\left( \prod_{a \in (\Z/p^n)^{\times}/\langle N \rangle } (1 - \zeta^a)^a
\right)^{2c}
= \left( \prod_{a \in (\Z/p^n)^{\times}/\langle N_1 \rangle}
(1 - \zeta^a)^a \right)^2$$
$$= \prod_{a \in (\Z/p^n)^{\times}/\langle N_1 \rangle}
(1 - \zeta^a)^a (1-\zeta^{-a})^{-a}
= \prod_{a \in (\Z/p^{\f1})^{\times}}(-\zeta)^{a^2}.$$
(Here the equalities hold in the quotient
$(\Of_n/\frak n)^{\times}/(p^n,N^2-1)$.)

If $p \geq 5,$ then since there are quadratic residues distinct from
1 in $(\Z/p)^{\times},$ we compute that
$\sum_{a \in (\Z/p^{\f1})^{\times}} a^2\equiv 0 \pmod{p^{\f1}},$
and so 
$\prod_{a \in (\Z/p^n)^{\times}/\langle N \rangle } (1 - \zeta^a)^a$
generates the trivial subgroup of 
$(\Of_n/\frak n)^{\times}/(p^n,N^2-1)$.  In this case, our ``refined analysis''
adds no further restrictions to the degree of $L_n$ over $K_n$.
However, if $p = 3$, then 1 is the only quadratic residue in
$(\Z/3)^{\times},$ and one computes that the power of $3$ dividing
$\sum_{a \in (\Z/p^{\f1})^{\times}} a^2$ is exactly
$3^{\f1 -1 } = 3^{f-1}$.  Thus we find that the degree 
$[L_n:K_n]$ is bounded above by $3^{f-1}.$  This is the exact power of
3 dividing $(N^2-1)/24$, and thus we have proved the lemma in the case
$p=3$.

Suppose now that $p=2$.   Write $N = \pm 1 \cdot N_1,$ where $N_1 \equiv 1
\pmod 4$.  Let $2^f$ be the exact power of 2 dividing $N^2-1,$
and let $2^{\f1}$ be the exact power of 2 dividing $N_1 - 1$.
Note that $f = \f1 + 1.$
Also, assume that $n \geq f$.  In particular, $n \geq 2$,
and so $-\zeta$ is also a primitive $2^n$th root of unity.
The quotient $(\Of_n/\frak n)^{\times}/(2^n,N^2-1)$ is then cyclic
of order $2^f$, generated by $\zeta,$ or by $-\zeta$.

We may rewrite 
$\prod_{a \in (\Z/2^n)^{\times}/\langle N \rangle } (1 - \zeta^a)^a$
in the form
$$\prod_{a \in (\Z/2^n)^{\times}/\langle N \rangle } (1 - \zeta^a)^a
=\prod_{a \in (1 + 4\Z/2^n)/\langle N_1 \rangle}
(1 - \zeta^a)^a (1-\zeta^{-a})^{-a}
= \prod_{a \in (1+4\Z/2^n)/(1+2^{\f1}\Z/2^n)} (-\zeta)^{a^2}.$$
One computes that the largest power of 2 dividing
$\sum_{a \in (1 + 4\Z/2^n)/(1+2^{\f1}\Z/2^n)} a^2$
is $2^{\f1-2} = 2^{f-3}$.  Thus the degree $[L_n:K_n]$ is bounded above
by $2^{f-3}$.
This is the exact power of $2$
dividing $(N^2-1)/24$, and so we have proved the lemma
in the case $p=2$.
$\qed$
\end{Proof}

\

{\em Conclusion of proof of Theorem~\ref{cohomology calculation}:}
Lemmas~\ref{one} and~\ref{two} yield a natural injection
$I/I^2 \iso \ilim A_n \hookrightarrow \ilim B_n.$
Together with Lemma~\ref{three}, 
this implies Theorem~\ref{cohomology calculation}.
$\qed$
\end{Proof}

\

The reduced Zariski tangent space
of the deformation ring $R$ can be computed via a calculation
similar to that used to prove Theorem~\ref{cohomology calculation}.
We state the result here, but postpone the details of the calculation
to the following sections.  (See Proposition~\ref{prop:tangent space at 2}
for the case $p=2$, and
Proposition~\ref{prop:tangent space at odd p}
for the case of odd $p$.)

\begin{prop}\label{tangent space}
If $\m$ denotes the maximal ideal of $R$, then
the reduced Zariski tangent space $\m/(\m^2,p)$ of $R$ is of
dimension at most one over $\F_p$.
More precisely,
$\m/(\m^2,p)$  vanishes unless
$p$ divides the numerator of $(N-1)/12$,
in which case it has dimension one over $\F_p$.
\end{prop}

Having introduced the deformation ring $R$, 
we now turn to constructing the corresponding Hecke ring $\T$.
We consider the space $M_2(N)$ of all modular forms of weight
two on $\Gamma_0(N)$ defined over $\Qbar_p$,
and the commutative $\Z_p$-algebra $H$ of endomorphisms
of $M_2(N)$ generated by the Hecke operators $T_n$.  
We define the $p$-Eisenstein maximal ideal of the algebra $H$ to be the
ideal generated by the elements $T_n - \sigma^*(n)$ (where
$
\sigma^*(n) = \sum_{{0 < d|n}\atop (d,N)=1} d$
for any positive integer $n$)
together with the prime $p$, and let $\T$ denote the completion
of~$H$ at its  $p$-Eisenstein maximal ideal.
Then $\T$ is a reduced $\Z_p$-algebra.
We let $J$ denote the kernel of the surjection $\T \rightarrow \Z_p$
describing the action of $\T$ on the Eisenstein series $E^*_2$, 
let $\T^0$ denote
the quotient of $\T$ that acts faithfully on cuspforms,
and let $J^0$ denote the image of $J$ in $\T^0$.
(This is the localisation at $p$ of the famous
Eisenstein ideal of \cite{eisenstein}.)

\begin{lemma}\label{congruence modulus}
The order of $\T^0/J^0$ $($which is a power of $p)$
is equal to the $p$-power part of the
numerator of $(N-1)/12$.
\end{lemma}

\begin{Proof} 
This is Proposition~II.9.7 of \cite{eisenstein}.
$\qed$
\end{Proof}

\begin{prop}\label{Hecke rep}
There is an object $(V,L,\rho)$ of $\Def(\T)$, uniquely
determined by the property that $\Trace(\rho(\Frob_{\ell}))
= T_{\ell}$, for $\ell \neq p,N$.  Furthermore,
the diagram
$$\xymatrix{R \ar[d] \ar[rr] && \T \ar[d] \\
R/I \ar@{=}[r] & \Z_p \ar@{=}[r] & \T/J}$$
is commutative.
\end{prop}

\begin{Proof}
Since, by Corollary~\ref{traces}, the universal deformation ring
$R$ is generated by the traces $\Trace(\rhouniv(\Frob_{\ell}))$,
there is at most one object $(V,L,\rho)$ of $\Def(\T)$ satisfying the
condition $\Trace(\rho(\Frob_{\ell})) = T_{\ell}$ for $\ell \neq p,N$.
This gives the uniqueness statement of the proposition.
In order to construct the required object $(V,L,\rho)$,
we proceed in several steps.  

\begin{lemma}\label{five}
Let $\Vbar'$ be a two dimensional continuous $G_{\Q}$-module
over a finite extension $k$ of $\F_p$.  Suppose that $\Vbar'$ is finite at $p$,
unramified away from $p$ and $N$, contains an $I_N$-fixed
line that is {\em not} $G_{\Q}$-stable, and has semi-simplification
isomorphic to the semi-simplification of $k\otimes_{\F_p}\Vbar.$
Then $\Vbar'\cong k\otimes_{\F_p} \Vbar.$
\end{lemma}

\begin{Proof}
Since $\Vbar'$ and  $k\otimes_{\F_p}\Vbar$ have isomorphic semi-simplifications,
we see that $\Vbar'$ is an extension of one of $k\otimes_{\F_p}\mu_p$ or
$k\otimes_{\F_p} \Z/p$
(thought of as \'etale groups schemes over $\Q$, or equivalently as 
$G_{\Q}$-representations) by the other.  Both these one dimensional
representations are unramified at $N$, and $\Vbar'$ contains
one or the other as a $G_{\Q}$-submodule.  It also contains 
an $I_N$-fixed line which is {\em not} a $G_{\Q}$-submodule.
Thus $\Vbar'$ is in fact spanned by $I_N$-fixed lines, and so
is unramified at $N$.  By assumption it is finite at $p$,
and so it has a prolongation to a finite flat group scheme over $\Z$. 

If $p$ is odd, then $\Vbar$ must prolong to an extension of one
of $k\otimes_{\F_p}\mu_p$ or $k\otimes_{\F_p} \Z/p$
by the other as a group scheme over
$\Spec \Z$  (since by \cite{fontaine}, Thm.~2, $p$-power order group
schemes over $\Z$ are determined by their associated Galois
representations).  There are no such non-trivial extensions
(\cite{eisenstein}, Ch.~I), and thus $\Vbar' \iso k\otimes_{\F_p}\Vbar.$
In the case that $p = 2,$ note first that since both
$k\otimes_{\F_2}\mu_2$ and $k\otimes_{\F_2}\Z/2$
yield the trivial character of $G_{\Q}$,
the module $\Vbar$ cannot be the direct sum of these
two characters; if it were, every line (including
the $I_N$-fixed line appearing in the statement of the lemma)
would be $G_{\Q}$-stable.  Taking this into account,
it is easily seen (again using the results of \cite{eisenstein},
Ch.~I) that $\Vbar' \iso k\otimes_{\F_2}\Vbar$.
$\qed$
\end{Proof}

\begin{lemma}\label{six}
Let $K$ be a finite extension of $\Q_p,$ with ring of integers $\Of$.
Let $k$ denote the residue field of $\Of,$ and let $\Of'$ denote
the order in $\Of$ consisting of elements whose image in $k$ 
lies in the prime subfield $\F_p$ of $k$.
Suppose given a two dimensional $K$-vector space $W$,
and a continuous representation  $G_{\Q} \rightarrow \GL(W)$ 
that is finite at $p$ $($in the sense that one, or equivalently
any, $G_{\Q}$-invariant $\Of$-lattice in $W$ is finite at $p)$,
semistable at $N$ $($in the sense that $W$ contains an $I_N$-fixed
line$)$, unramified away from $p$ and $N$, such that the semi-simple 
residual representation attached to $W$ is isomorphic to the
direct sum of the trivial character and the mod $p$ cyclotomic character.

If $W$ is irreducible, then we may find a free $\Of'$-module of rank
two $V$, equipped with a continuous representation $\rho:
G_{\Q} \rightarrow \GL(V)$, and containing an $I_N$-fixed line $L$,
such that the triple $(V,L,\rho)$ deforms $(\Vbar,\Lbar,\rhobar)$,
and such that $K\otimes_{\Of'} V \iso W$ as $G_{\Q}$-modules.
\end{lemma}

\begin{Proof}
Choose any $G_{\Q}$ lattice $V'$ in $W$, and let $L'$ denote
the intersection of $V'$ with the $I_N$-fixed line in $W$.
Since $W$ is irreducible, the line $L'$ is not $G_{\Q}$-stable.
Thus we may find a natural number $n$ such that
$L'/p^n$ is $G_{\Q}$-stable in $V'/p^n,$ but such that
$L'/p^{n+1}$ is {\em not} $G_{\Q}$-stable in $V'/p^{n+1}$.
If we define $V''$ to be the preimage in $V'$ of $L'/p^n$,
then we see that $L'/p$, when regarded as a subspace of
$V''/p,$ is not $G_{\Q}$-stable.  
Lemma~\ref{five} implies that $V''/p \cong k\otimes_{\F_p}\Vbar$.
Using the description of the automorphisms of $k\otimes_{\F_p}\Vbar$
afforded by Lemma~\ref{minendos}, we deduce easily that in
fact there is an isomorphism of pairs $(V''/p,L'/p) \iso k\otimes_{\F_p}
(\Vbar,\Lbar).$  If we choose a basis for $(V'',L')$ over $\Of$ that reduces
to an $\F_p$ basis for $(\Vbar,\Lbar)$, then the $\Of'$-span of
this basis gives rise to the required pair $(V,L)$.
$\qed$
\end{Proof}

\

If $\twT$ denotes the normalisation of $\T$, then we may
write $\twT = \prod_{i=1}^d \Of_i,$
where each $\Of_i$ is a discrete valuation ring, of finite
index over $\Z_p$.   The rings $\Of_i$ are in bijection
with the conjugacy classes of normalised eigenforms $f_i$ in
$M_2(N)$ that satisfy the congruence $f_i \equiv E^*_2 \pmod \m_i$
(where $\m_i$ denotes the maximal ideal of $\Of_i$); as before,
$E^*_2$ denotes the weight two Eisenstein series on $\Gamma_0(N)$.
The ring $\Of_i$ is the ring of integers in the subfield of $\Qbar_p$
generated by the Fourier coefficients of $f_i$.
The injection $\T \rightarrow \twT = \prod_{i=1}^d \Of_i$ is characterised
by the property
$T_n \mapsto (a_n(f_i))_{i=1,\ldots,d}.$
Note that $E^*_2$ is one such form $f_i$.  We may choose the
labeling so that $E^*_2 = f_1$; then $\Of_1 = \Z_p = \T/J$.

As in the statement of Lemma~\ref{six}, for each $i = 1,\ldots,d$,
define $\Of_i'$ to be the order in $\Of_i$
obtained as the preimage under the map to the residue field 
of the prime subfield $\F_p$.
By construction $\Of_i'$ is a complete Noetherian local
ring with residue field $\F_p$.  Also, the natural map $\T \rightarrow \Of_i$
factors through $\Of_i'$.

\begin{lemma}\label{seven} For each $i = 1,\ldots,d$,
we may construct an object $(V_i,L_i,\rho_i) \in \Def(\Of_i')$
with the property that $\Trace(\rho_i(\Frob_{\ell}))$ is equal
to the image of $T_{\ell}$ in $\Of_i'$, for each $\ell \neq p,N$.
\end{lemma}

\begin{Proof}  If $i = 1$, so that $\Of_i' = \Z_p,$ we take
$(V_1,L_1,\rho_1)$ to be the triple $(\Vmin,\Lmin,\rhomin).$
Suppose now that $i\geq 2,$ so that $\Of_i$ corresponds to a cuspform
$f_i$.  If we consider the usual irreducible Galois representation
into $\GL_2(\Q_p\otimes_{\Z_p} \Of_i)$ attached to $f_i$,
and apply Lemma~\ref{six}, then we again obtain the required triple.
$\qed$
\end{Proof}

\

{\em Conclusion of proof of Proposition~\ref{Hecke rep}:}
Each of the triples $(V_i,L_i,\rho_i)$ constructed in the
previous lemma corresponds to a homomorphism $\phi_i : R \rightarrow \Of_i'.$
The product of all these yields a homomorphism
$\phi: R \rightarrow \prod_{i = 1}^d \Of_i'.$ 
Since $R$ is topologically generated by the elements 
$\Trace(\rhouniv(\Frob_{\ell}))$ ($\ell \neq p,N$),
we see that $\phi$ factors through $\T$.  The map $\phi$ in turn
corresponds to a triple $(V,L,\rho) \in \Def(\T),$
satisfying the requirements of the proposition.
By construction, the diagram appearing
in the statement of the proposition commutes.
$\qed$
\end{Proof}

\

Let $\T'$ denote the image in $\T$ of the map constructed
in Proposition~\ref{Hecke rep}.  Our ultimate goal
is to prove that that map is an isomorphism, and so in
particular that $\T' = \T$.  However, we will proceed in
stages.

Write $J' = \T' \bigcap J$,
let $(\T')^0$ denote the image of $\T'$ in $\T^0$,
and let $(J')^0$ denote the image of $J'$ in $(\T')^0$.
We have the morphism of short exact sequences
$$\xymatrix{ 0 \ar[r] &\T' \ar[r]\ar[d] & \Z_p \bigoplus (\T')^0 \ar[r]\ar[d] &
(\T')^0/(J')^0 \ar[r]\ar[d] & 0 \\
0 \ar[r] & \T \ar[r] & \Z_p \bigoplus \T^0 \ar[r] &
\T^0/J^0 \ar[r] & 0 .}$$
Applying the snake lemma we obtain
 the following exact sequence:
$$
\begin{matrix} 0 \longrightarrow \kernel((\T')^0/(J')^0 \rightarrow \T^0/J^0)
\longrightarrow \cokernel(\T' \rightarrow \T) \qquad \qquad \qquad
& \\ \qquad \qquad \qquad
\longrightarrow \cokernel((\T')^0 \rightarrow \T^0)
\longrightarrow
\cokernel((\T')^0/(J')^0 \rightarrow \T^0/J^0) \longrightarrow 0.\end{matrix}
$$
We also have the following tautological exact sequence:
$$\begin{matrix} 0 \longrightarrow\kernel((\T')^0/(J')^0 \rightarrow \T^0/J^0)
\longrightarrow (\T')^0/(J')^0 \qquad
& \\ \qquad \qquad \qquad
\longrightarrow \T^0/J^0
\longrightarrow
\cokernel((\T')^0/(J')^0 \rightarrow \T^0/J^0)\longrightarrow 0.\end{matrix}$$
Thus we find that 
\begin{equation}\label{inequality}
\#(\T')^0/(J')^0 - \#(\T^0/J^0) = \# \cokernel(\T' \rightarrow \T) -
\# \cokernel((\T')^0 \rightarrow \T^0).
\end{equation}
Since $\T \rightarrow \T^0$ is surjective, we conclude that
the right hand side of~(\ref{inequality}) is non-negative,
and thus that the order of $(\T')^0/(J')^0$ is at least equal to
that of $\T^0/J^0$.  By Lemma~\ref{congruence modulus},
the order of this latter group has order equal to the $p$-power part
of the numerator of $(N-1)/12.$   Thus the order of
$(\T')^0/(J')^0$ is at least equal to this number.

Suppose now that $N \not\equiv -1 \pmod{2p}$.  
The $p$-power part of $(N^2-1)/24$ is then equal
to the $p$-power part of the numerator of $(N-1)/12$.
The numerical criterion of \cite{wiles}
(as strengthened in \cite{lenstra}) thus applies
to show that the surjection $R \rightarrow \T'$ of the
preceding proposition is an isomorphism of local complete
intersections. 
Furthermore, we conclude that in fact $(\T')^0/(J')^0$
has order exactly equal to the power of $p$ dividing the
numerator of $(N-1)/12$, that is, to the order of $\T^0/J^0$.
The equation~(\ref{inequality}) then shows
that $\T' = \T$ if and 
only if $(\T')^0 = \T^0.$

\begin{lemma}\label{eight} The inclusion $\T' \rightarrow \T$
is an isomorphism.
\end{lemma}

\begin{Proof} It follows from Corollary~\ref{traces}, together with
the construction of the map $R \rightarrow \T$ of
Proposition~\ref{Hecke rep}, that $\T'$ contains $T_{\ell}$
for all $\ell \neq N, p$.  Proposition~\ref{ordinary} shows
that $\rhouniv$ has a rank one space of $I_p$-coinvariants,
on which $\Frob_p$ then acts as multiplication by a unit.  It follows from
the construction of $R \rightarrow \T$, and the known structure
of Galois representations attached to modular forms, that
the image of this unit in $\T$ is equal to the Hecke operator~$T_p.$ 
Thus $\T'$ contains $T_p$.

It remains to show that $T_N$ lies in $\T'$.  By the remark
preceding the statement of the lemma, it in fact
suffices to show that $T_N$ lies in $(\T')^0$.   
The surjection $R \rightarrow (\T')^0$ induces an object
$(V^0,L^0,\rho^0) \in \Def((\T')^0)$.
The concrete construction of the map $R \rightarrow \T$
(and hence the map $R \rightarrow \T^0$) shows that 
this representation is built out of Galois representations
attached to cuspforms on $\Gamma_0(N)$, which are
(so to speak) genuinely semi-stable at $N$.  In particular,
the line $L^0$ is not only fixed by 
$I_N$, but is stable under the decomposition group at $N$.
Standard properties of Galois representations attached to cusp forms
show that the eigenvalue of $\Frob_N$ on this line is furthermore
equal to $T_N$.  Thus $T_N \in (\T')^0,$ and so we see
that $(\T')^0 = \T^0,$ as required.
$\qed$
\end{Proof}

\

The preceding lemma completes the proof of Theorem~\ref{R=T}
in the case when $N \not \equiv -1 \pmod{2p}$.
If, on the other hand, we have $N \equiv -1 \pmod{2p},$
then Proposition~\ref{tangent space} shows that
the Zariski tangent space of $R$ is trivial.
In this case, the map $R\rightarrow \Z_p$ is an isomorphism.
Also, Lemma~\ref{congruence modulus} then implies that
$\T^0 = 0,$ and hence that $\T = \Z_p$.  Thus the
map $R \rightarrow \T$ is certainly an isomorphism in this case,
and we have completely proved Theorem~\ref{R=T} of the introduction.

\

Let us make two remarks:

\smallskip

(A) An alternative approach to proving Proposition~\ref{Hecke rep} is
as follows.  The results of \cite{eisenstein}, Section~II.16,
show that if $V^0$ denotes the $p$-Eisenstein part of the $p$-adic Tate module
of $J_0(N)$, then $V^0$ is free of rank two over $\T^0$, and 
the $G_{\Q}$-action on $V^0$ 
yields a deformation $\rho^0$ of $\rhobar$ over $\T^0$.  The $I_N$-invariants
in this representation form a rank one free submodule $L^0$ of this
representation.   The discussion of \cite{eisenstein}, Section~II.11
shows that both the cuspidal and Shimura subgroup map isomorphically
onto the connected component group of the fibre over $N$ of
the N\'eron model of $J_0(N)$, and this in turn implies that
$(V^0,L^0,\rho^0)$ provides an object of $\Def(\T^0)$.
Thus we obtain a corresponding map $R \rightarrow \T^0$.  Taking the
product of this with
the map $R \rightarrow R/I = \Z_p,$ we obtain the required map
$R \rightarrow \T$ of Proposition~\ref{Hecke rep}.
Finally, the explicit description of $\T^0$ provided by
\cite{eisenstein}, Cor.~II.16.2 assures us that the
map $R \rightarrow \T$ is surjective.

We have chosen to present the alternative argument above
both because it is more elementary (the only ingredient
required from \cite{eisenstein}, Ch.~II, is the
computation of the order of $\T^0/J^0$), and because
we are then able to recover the results of
\cite{eisenstein}, Sections~II.16, II.17, as we
explain below.

\smallskip

(B) In the proof of Lemma~\ref{eight}, we have struggled
slightly to prove that $T_N$ in fact lies in $\T'$.
This is somewhat amusing, since in fact 
$T_N = 1$ in $\T$!  This follows from \cite{eisenstein},
Prop.~II.17.10.  When $p$ is odd, the argument is 
straightforward: namely, since $T_N^2 = 1$ for general
reasons (the Galois representations attached to
modular forms on $\Gamma_0(N)$ are semi-stable at 
$N$ and Cartier self-dual), it suffices to note
that $T_N \equiv 1$ modulo the maximal ideal of $\T$.
When $p = 2$, Mazur's proof of this result depends
on his detailed analysis of the 2-Eisenstein torsion in 
$J_0(N)$.  We present an alternative proof below,
using the deformation theoretic techniques of
this paper.

\

We close this section by explaining how Theorem~\ref{R=T}
allows us to recover the main results of Section II of~\cite{eisenstein}.

\begin{cor}\label{gorenstein} The $\Z_p$-algebra $\T$
$($and consequently also its
quotient $\T^0)$ is generated by a {\em single} element over $\Z_p$.
In particular, both $\T$ and $\T^0$ are local complete intersections,
and hence Gorenstein.
\end{cor}

\begin{Proof} 
Theorem~\ref{R=T} shows
that it suffices to verify the analogous statement for the
deformation ring $R$.  Proposition~\ref{tangent space} shows
that if $\m$ denotes the maximal ideal of $R$,
then $\m/(\m^2,p)$ has dimension at most one over $\F_p$,
and the corollary follows.
$\qed$
\end{Proof}

\

The fact that $\T^0$ is monogenic over $\Z_p$ was originally proved by Mazur
(\cite{eisenstein}, Cor.~16.2).  
Since $\T$ is monogenic over $\Z_p$, and is equipped with a map
$\T \rightarrow \T/J \iso \Z_p,$ we see that we may write 
$\T \iso \Z_p[X]/Xf(X),$ where $X$ generates the ideal $J$ in $\T$,
the polynomial $f(X) \in \Z_p[X]$ satisfies
$f(X) \equiv X^{g_p} \pmod p$, and there is an isomorphism
$\T^0 \cong \Z_p[X]/f(X).$
(Here we follow \cite{eisenstein} in letting $g_p$ denote the rank
of $\T^0$ over $\Z_p$.)  The image of $X$ in $\Z_p[X]/f(X)$
generates the ideal $J^0$ in $\T^0$.  

In \cite{eisenstein}, Prop.~II.18.10, Mazur treats
the questions of exhibiting explicit generators of $J^0$
(or equivalently, explicit choices for the element
``$X$'' of the preceding paragraph).  We recall his
result here, and give a deformation-theoretic proof.

\begin{prop}\label{prop:good primes}
Suppose that $p$ divides the numerator of $(N-1)/12$.
Let $\ell$ be a prime different from $N$.
Say that $\ell$ is {\em good} $($with respect to
the pair $(p,N))$ if (i) one of $\ell$ or $p$
is odd, $\ell$ is not a $p$th power modulo
$N$, and $(\ell - 1)/2 \not\equiv 0 \pmod p;$
or (ii) $\ell = p = 2$ and $-4$ is not
an $8$th power modulo $N$.
\footnote{This definition originally
appeared in \cite{eisenstein}, p.~124.  However,
condition~(ii) is misstated there.}

Then $T_{\ell} - (1 + \ell)$ generates the
ideal $J^0$ if and only if $\ell$ is a good prime.
\end{prop}

\begin{Proof}
Let $R \cong \T \rightarrow \F_p[X]/X^2$ be a map that classifies
a (unique up to scaling, 
by Proposition~\ref{tangent space}) non-trivial element
in the reduced Zariski tangent space of $R$.
If $\ell$ is distinct from $p$, then
we must show that $T_{\ell} - (1 + \ell)
= \Trace(\rhouniv(\Frob_{\ell})) - (1 + \ell)$ has non-zero
image under this map if and only if $\ell$ is a good prime.
If $\alpha_p\in R \cong \T$ denotes the
scalar by which $\Frob_p$ acts on the rank one quotient module
of $I_p$-coinvariants of $\Vuniv$, then $T_p = \alpha_p$, and 
so we must also show that $\alpha_p - (1 + p)$ has non-zero image
under this map if and only if $p$ is a good prime.
Both cases follow from Proposition~\ref{prop:tangent space at 2}
in the case when $p=2$, and from
Proposition~\ref{prop:tangent space at odd p}
in the case of odd $p$.
$\qed$
\end{Proof}

\


As was remarked upon above, the next result (and the final result
of this section) is also originally due to Mazur.

\begin{prop}\label{T_N = 1}
In $\T$ we have the equality $T_N = 1$.
\end{prop}

\begin{Proof} As we recalled above, this result is straightforward
when $p$ is odd.  Thus we assume that $p = 2$.  The $T_N$-eigenvalue
of $E^*_2$ is equal to 1.  Thus, in order to
show that $T_N = 1$, it suffices to show that for each
cuspform $f_i$ ($i = 2,\ldots,d$ -- we are using the notation
introduced during the proof of Proposition~\ref{Hecke rep}),
the image of $T_N$ in $\Of_i$ is equal to $1$.
If $N \not\equiv 1 \pmod 8$, then there are no cuspforms to
consider, and hence there is nothing to prove.  Thus we assume
for the remainder of the argument that $N \equiv 1 \pmod 8$.

Fix a cuspform $f_i$, and let $S$ denote the local ring
$$S = \{(a,b) \in \Z/4 \times \Of_i/2\m_i \, | \, a  \bmod 2 = b \bmod \m_i\}.$$
The objects $(\Vmin_2,\Lmin_2,\rhomin_2) \in \Def(\Z/4)$ and
the object in $\Def(\Of_i'/2\m_i)$ obtained by reducing modulo
$2\m_i$ the object
$(V_i,L_i,\rho_i) \in \Def(\Of_i')$ (the latter was constructed in the
course of proving Proposition~\ref{Hecke rep}) glue to yield
an object $(V,L,\rho) \in \Def(S)$.
Since $N \equiv 1 \pmod 8,$ we see that $G_{\Q_N}$ acts trivially
on $\Vmin_2$.  Since $(V_i,L_i,\rho_i)$ is constructed from
the Galois representation attached to the cuspform $f_i,$ we
see that $G_{\Q_N}$ stabilises $L_i,$ and $\Frob_N$ acts
as multiplication by $T_N$ on $L_i$.  Thus the line $L$ is stabilised by $G_{\Q_N}$
(in addition to being fixed by $I_N$), and $\Frob_N$ acts
as multiplication by the image of $T_N$ in $S$.  If the image of $T_N$ in 
$\Of_i$ is equal to $-1$, then we see that the image of $T_N$ 
in $S$ is equal to $(1,-1)$.
Now $(1,-1) \not\equiv (1,1) \pmod{2S}$.
Thus the object $(V/2,L/2,\rho/2) \in \Def(S/2)$ obtained by
reducing $(V,L,\rho)$ modulo 2 has the property that
$L$ is stable, but not trivial, under the action of $G_{\Q_N}$.
On the other hand, Theorem~\ref{p=2 class field
factorisation}, together with Lemma~\ref{lem:info at N}, shows
that there are no elements of $\Def(S/2)$.  This contradiction
proves the proposition.
$\qed$
\end{Proof}

\section{Explicit deformation theory: $p$ = 2}
\label{sec:explicit at 2}
Let us begin by fixing an odd prime $N$,
and recalling some class field theory
of the field $K = \Q(\sqrt{(-1)^{(N+1)/2}N}).$
We let $H$ denote the 2-power part of the
strict class group $\Cl(\Of_K)$ of the ring of integers $\Of_K$ of $K$,
and let $E$ denote the corresponding cyclic 2-power
extension of $K$, which is unramified at all
finite primes.
Genus theory shows that $H$ is cyclic,
and non-trivial.  Thus $E$ is a non-trivial
cyclic 2-power extension of $K$; its unique
quadratic subextension is equal to $K(\sqrt{-1}).$
We let $\pi_K$ denote the unique prime
of $K$ lying over 2; its image in $H$
generates the two-torsion subgroup $H[2]$ of $H$.

The following result is classical, but we will
recall a proof for the benefit of the reader.

\begin{prop}\label{prop:2-cft}
The order of $H$ is divisible by four if and only
if $N \equiv 1 \pmod 8.$  The order of $H$ is 
divisible by eight if and only if furthermore
$-4$ is an $8$th power modulo $N$.
\end{prop}

\begin{Proof}
If $N \equiv -1 \pmod 4,$ then $K$ is a real
quadratic field.  If $E^+$ denotes the 2-Hilbert
class field of $K$ (so $E^+$ is the maximal
totally real subextension of $E$), then
we see that $E$ is equal to the compositum
of $E^+$ and $K(\sqrt{-1}).$  Since $E$ is
cyclic over $K$, we deduce that $E^+$ must
in fact be trivial.  Thus in this case
$H$ is of order two.

Suppose now that $N \equiv 1 \pmod 4,$ and
that $E$ contains a degree four sub-extension.
Since $E/K$ is cyclic, this sub-extension is
unique, and hence Galois over $\Q$.
It must contain $\Q(\sqrt{-1})$,
and one sees easily that it is in fact 
a biquadratic extension of $\Q(\sqrt{-1})$,
unramified away from $N$.
Since it is Galois over $\Q$, it must be
of the form
$\Q(\sqrt{-1},\sqrt{\nu}, \sqrt{\nubar}),$
where $\nu$ is an element of $\Z[\sqrt{-1}]$
(and $\nubar$ is its conjugate) satisfying
$\nu\nubar = N.$

However, for the extension
$\Q(\sqrt{-1},\sqrt{\nu},\sqrt{\nubar})/\Q(\sqrt{-1})$
to actually be unramified at 2, it must be that
$\nu \equiv 1 \pmod 4.$  The element $\nu$ can
be chosen in this manner if and only
if $N \equiv 1 \pmod 8.$  Thus we see that $E$
has a degree four subfield if and only
if $N$ satisfies this congruence.

Finally, let us consider the question of whether
the order of $H$ is divisible by eight.  This
is the case if and only if the two-torsion
subgroup $H[2]$ of $H$ has trivial image in
$H/4$; equivalently, if and only if $\pi_K$
has trivial image in $H/4$.
This holds, in turn,
if and only if $\pi_K$ splits
completely in 
$\Q(\sqrt{-1},\sqrt{\nu}, \sqrt{\nubar}).$
Clearly, this is true if and only if 
the ideal $(1+i)$ splits completely
in this field, regarded as an extension of $\Q(\sqrt{-1}).$
This holds, in turn, if and only if
$(1+i)$ is a quadratic residue modulo $\nu$
(or equivalently module $\nubar$).  
Raising to $4$th powers, and taking
into account the isomorphism
$\Z/N \iso \Z[\sqrt{-1}]/\nu,$
we see that this is equivalent
to $-4$ being an $8$th power modulo $N$.
$\qed$
\end{Proof}

\

The following lemma is used in
the proof of Proposition~\ref{T_N = 1}.

\begin{lemma}\label{lem:info at N}
The inertia group $I_N$ and the decomposition group $G_{\Q_N}$
have the same image in $\Gal(E/\Q).$
\end{lemma}

\begin{Proof}
There is a unique prime lying above $N$ in $K$,
and it principal.  Thus this prime splits completely in the
Hilbert class field $E$ of $K$.  
Thus $I_N$ and $G_{\Q_N}$ both have trivial image in $\Gal(E/K)$.
Since $N$ is ramified in $K/\Q$, the lemma follows.
$\qed$
\end{Proof}

\

Let $H'$ denote the $2$-power part of the
strict ray class group of $K$ of conductor $\pi_K^2,$
and let $H''$ denote the $2$-power part of the
strict ray class group of $K$ of conductor $\pi_K^3$.
(Here ``strict'' means that in the case when $K$ is real quadratic, we
allow ramification at infinity.)  We let $E'$ and $E''$ denote
the corresponding abelian extensions of $K$.

\begin{prop}\label{prop:RCF} (i) The natural surjection $H'' \rightarrow H'$
is an isomorphism. 

(ii) Either the natural surjection $H' \rightarrow H$ is an isomorphism,
in which case $E = E'$,
or else the kernel of this surjection has order two,
in which case $E'/E$ is a quadratic extension that is ramified at two.

(iii) The group $H'$ is cyclic. 

(iv) 
Let $D_2(E'/\Q)$ denote the decomposition group of some prime
of $E'$ lying over $2$, and let $I_2(E'/\Q)$ denote the
inertia subgroup of $D_2(E'/\Q)$. 
Then $I_2(E'/\Q)$ has index two in $D_2(E'/\Q)$.
If furthermore $E'/E$ is a quadratic extension,
then $D_2(E'/\Q)$ is dihedral of order $8$.
\end{prop}

\begin{Proof}
The groups $H'$ and $H''$ sit inside the
following exact diagram:
$$\xymatrix{
\Of^{\times}_K \ar[r]\ar@{=}[d] & (\Of_K/\pi^3_K)^{\times} \ar[r]\ar[d]^{\psi}
& H'' \ar[r] \ar[d] & H \ar[r] \ar@{=}[d]& 0 \\
\Of^{\times}_K \ar[r] & (\Of_K/\pi^2_K)^{\times}
\ar[r] & H' \ar[r] 
& H \ar[r] & 0.}$$ 
To prove that the map $H''\rightarrow H'$ is an isomorphism,
it suffices to show that the kernel
of $\psi$ maps to zero in $H''$; in other
words, that the kernel of $\psi$ consists of the images of global units.
Since $\pi^2_K = (2)$, we see that the kernel of $\psi$ is
equal to $\{\pm 1\}$; this completes the proof of~(i).

The proof of~(ii) is even more straightforward: it follows
immediately from a consideration of the exact sequence
$$\Of_K^{\times} \rightarrow (\Of_K/\pi_K^2)^{\times}
= (\Of_K/2)^{\times} \rightarrow H' \rightarrow H \rightarrow 0.$$

We now turn to proving~(iii).  For this, it
suffices to prove that $H'/2  \cong \Z/2 $.
Note that since the non-trivial element of $\Gal(K/\Q)$ acts on $H'$ 
via $h \mapsto h^{-1}$, we see that the extension $K'$ of $K$
corresponding by class field theory to $H'/2$ is abelian
over $\Q$.  If $H'/2$ were isomorphic to $\Z/2 \oplus \Z/2$
(rather than $\Z/2$),
then since $H$ is cyclic,
this would imply that there
exists a subextension of $K'$, quadratic over $K$ and of conductor $2$.
Such an extension would again be abelian over $\Q$.
Using the Kronecker-Weber theorem,
it is easy to check that there are no quadratic extensions of $K$.
Thus $H'$ must be cyclic, as claimed.

As remarked upon above, the class of $\pi_K$ has order two
in $H$.
Thus the decomposition group $D_2(E/K)$ at $2$ in the Hilbert
class field has order exactly $2$.  Since $K/\Q$ is ramified
at $2$, we see that the decomposition group $D_2(E/\Q)$ has
order four, and that the inertia subgroup $I_2(E/\Q)$ has order two.
If $E' = E$ then this completes the proof of~(iv).
If instead $E'/E$ is quadratic, then $E'/E$ is ramified at $2$,
implying that $D_2(E'/\Q)$ has order $8$, and that $I_2(E'/\Q)$
has order $4$.
Since $D_2(E'/K) \subseteq \Gal(E'/K) \cong H'$ is cyclic, by~(iii),
and since $\Gal(E'/\Q)$ is dihedral, it follows that $D_2(E'/\Q)$
is dihedral of order $8$.  
$\qed$
\end{Proof}

\

Let $(V,L,\rho)$ be an object of $\Def(A)$ for some
Artinian local $\F_2$-algebra, and let $F$ denote compositum
of $K$ with the fixed field of the kernel of $\rho$.  The following
result greatly restricts the possibilities for~$F$.

\begin{theorem}\label{p=2 class field factorisation}
The field $F$ is contained in the
strict $2$-Hilbert class field $E$ of $K$.
\end{theorem}

\begin{Proof}
Since $A$ is assumed to be of characteristic $2$,
the natural map $\Z_2^{\times} \rightarrow A^{\times}$
has trivial image, and thus the image of $\rho$ is contained in
$\SL_2(A)$. 
Since $I_N$ acts trivially on both $L$ and $V/L$,
we deduce that inertia at $N$
acts through an abelian group of exponent $2$, and thus 
through a cyclic group of order at most $2$.

\begin{lemma} 
The extension $F/K$ is unramified at all finite primes outside $\pi_K$.
Moreover, if $\Kab/K$ is the maximal abelian extension
of $K$ contained in $F$, then
the finite part of the conductor of $\Kab/K$ divides  $\pi^3_K$.
\end{lemma}

\begin{Proof}  The Galois group $\Gal(F/\Q)$ 
embeds into $\Gal(K/\Q) \times \rho(G_{\Q})$.  In particular, it is of
2-power order, and so
the image of an inertia group $I_N$ at $N$ in $\Gal(F/\Q)$ is cyclic
of two-power order.  As observed above, $\rho(I_N)$ is a quotient of
$I_N$ of order at most two.  On the other hand,
since $K/\Q$ is a quadratic extension that is ramified at $N$,
we see that $I_N$ surjects onto the order two group $\Gal(K/\Q)$.
It follows that the image of $I_N$ in $\Gal(F/\Q)$ has trivial
intersection with $\Gal(F/K)$, and so $F/K$ is unramified
at the prime above $N$.

By definition, $\rho$ is unramified outside $2$ and $N$,
and so it remains to prove the result about the conductor
of $\Kab/K$. Since the compositum of extensions of
conductor dividing $\pi^3_K$ has conductor at most $\pi^3_K$, it
suffices to prove the result for extensions
$K'/K$ with \emph{cyclic} Galois group. Suppose
such an extension $K'/K$ with Galois group $\Z/2^k \Z$
had conductor divisible by $\pi^4_K$. Then  the
conductor discriminant formula says
that the discriminant $\Delta_{K'/K}$ is the
product over all characters of $\Z/2^k \Z$ of
the corresponding conductor:
$$\Delta_{K'/K} = \prod_{\chi} \frak{f}_{\chi}.$$
 Since $\Z/2^k \Z$  has exactly
$2^{k-1}$ faithful characters, restricting the
product to this set we find that the discriminant
is divisible by at least
$(\pi_K)^{4 \cdot 2^{k-1}}$. In particular, this 
implies a lower bound for the two part of the
root discriminant of $K'$, and thus of $F$.
Explicitly,
$$\delta_{2,F} \ge \delta_{2,K'} = \delta_{2,K} N_{K/\Q}(\Delta_{K'/K})^{1/[K':\Q]}
\ge 2 \cdot 2 = 4.$$
Yet the Fontaine bound (\cite{fbound}, Theorem 1)
 for finite flat group schemes
over $\Z_2$ killed by $2$ implies that $\delta_{2,F}
< 2^{1 + \frac{1}{2-1}} = 4$. Thus the result follows
by contradiction. 
$\qed$
\end{Proof}

\

We will strengthen this lemma step-by-step,
until we eventually establish the theorem.

\begin{lemma}\label{lem:intermediate} The extension $F/K$ is cyclic, and
is contained in the field $E'$.
\end{lemma}

\begin{Proof}
The preceding lemma, together with part~(i) of
Proposition~\ref{prop:RCF},
shows that the extension of $K$ cut out by
any abelian quotient
of $\Gal(F/K)$ is contained in $E'' = E'$.
Part~(iii) of the same proposition then
implies that any such quotient is cyclic.
Thus $\Gal(F/K)$ is a $2$-group with no
non-cyclic abelian quotients, and so is itself
cyclic.  The result follows.
$\qed$
\end{Proof}

\


We now turn to a more careful study of the ramification
at $2$. Corollary~\ref{uniquenesscor} shows that
$V_{/\Q_2}$ has a unique prolongation to
a finite flat group scheme $M_{/\Z_2}$, that the action of $A$ on $V$
prolongs to an action of $A$ on $M$, and that
the connected-\'etale sequence
$$0 \rightarrow M^{\co} \rightarrow M \rightarrow M^{\et}
\rightarrow 0$$
induces a two-step filtration of $V$ by free $A$-modules
of rank one.

\begin{lemma} The action of inertia at $2$ on $ M^{\co}(\Qbar_2)$
and $ M^{\et}(\Qbar_2)$ is trivial.
\end{lemma}

\begin{Proof} This is clear for  $M^{\et}(\Qbar_2)$,
since \'{e}tale implies unramified. It follows for
$ M^{\co}(\Qbar_2)$ from the Cartier self-duality of $M_{/\Qbar_2}$.
$\qed$
\end{Proof}

\begin{lemma} If $\sigma \in G_{\Q_2}$ then
$\sigma^2$ acts trivially  on $V$.
\end{lemma}

\begin{Proof} Let us choose a basis of $V$ compatible with its
filtration arising from the connected-\'etale sequence
of $M$, and write the action of $\sigma$ on
$V$ as a matrix over $A$ in terms of this basis:
$$\sigma = \left(\begin{array}{cc} 1 + a & b  \\
0 & 1 + c \end{array} \right)$$
Part~(iv) of Proposition~\ref{prop:RCF} implies
that $\sigma^2$ lies in the inertia subgroup.
Thus it must fix $ M^{\co}(\Qbar_2)$
and $ M^{\et}(\Qbar_2)$. Computing $\sigma^2$,
we find that $(1+a)^2 = (1+c)^2 = 1$, and so $a^2 = c^2 = 0$.
Since the determinant of $\sigma$ is $1$, we see that
$(1+c) = (1+a)^{-1} = 1 - a$. Now computing $\sigma^2$ we
find that it is trivial.
$\qed$
\end{Proof}

\

{\em Conclusion of proof of Theorem~\ref{p=2 class field factorisation}:}
If $E' = E$,
then by Lemma~\ref{lem:intermediate} there is nothing more to prove.
Otherwise, Proposition~\ref{prop:RCF} implies that
the $D_2(E'/\Q)$ is dihedral of order $8$.
We have seen that for any $\sigma \in G_{\Q_2}$, the element
$\sigma^2$ acts trivially. Thus the image $\rho_{|G_{\Q_2}}$ factors through
an exponent $2$ group, which is therefore abelian. Yet
the dihedral group of order $8$ is not abelian, and hence
$F$ is contained in a proper subfield of $E'$ that is Galois
over $\Q$. All such subfields lie inside $E$.
$\qed$
\end{Proof}

\

\begin{cor}\label{cor:existence}  If $2^m$ denotes the order of $H$, then 
there exists a surjection $R \rightarrow \F_2[X]/X^n$
if and only if $n \leq 2^{m-1}$.  Furthermore, any
such surjection is unique up to applying an automorphism
of $\F_2[X]/X^n$.
\end{cor}

\begin{Proof}
Corollary~\ref{traces} implies that
there exists a surjection $R \rightarrow \F_2[X]/X^n$ if and only
if there exits $(V,L,\rho) \in \Def(\F_2[X]/X^n)$ with
the property that the traces of $\rho$ generate $\F_2[X]/X^n$
as an $\F_2$-algebra
(or equivalently, with the property that there is 
an element of $G_{\Q}$ whose image under $\rho$ has
trace congruent to $X \pmod{X^2}$.)

\begin{lemma}\label{lem:2-powers}
Let $A$ be an $\F_2$-algebra, and let $U \in \SL_2(A)$.

(i) $U^2 = I + \Trace(U) U.$

(ii) For any $k \geq 1$, we have that $\Trace(U^k) \in \Trace(U) A.$

(iii)
If $U \in \SL_2(A),$ then
$$U^{2^k} = \left( \sum_{i=0}^{m-1} \Trace(U)^{2^k-2^{k-i}} \right) I
+ \Trace(U)^{2^k-1} U,$$
for any $k \geq 1$.
\end{lemma}

\begin{Proof}
Any $2\times 2$ matrix $U$ over the ring $A$
satisfies the identity $U^2  = \Det(U) I + \Trace(U) U.$
Part~(i) is a particular case of this identity, and parts~(ii) and~(iii)
follow by induction.
$\qed$
\end{Proof}

\

Theorem~\ref{p=2 class field factorisation}
shows that $\rho$ factors as
$G_{\Q} \rightarrow \Gal(E/\Q) \rightarrow \SL_2(\F_2[X]/X^n).$
Now $\Gal(E/\Q)$ is a dihedral group of order $2^{m+1}$;
indeed, we may write
\begin{equation}\label{eq:presentation}
\Gal(E/\Q)= \langle \sigma, \tau | \sigma^{2^m} = \tau^2 =
(\sigma \tau)^2 = 1 \rangle,
\end{equation}
where $\sigma$ generates $\Gal(E/K)$, and $\tau$ generates the image of $I_N$
in $\Gal(E/\Q).$

Part~(i) of Lemma~\ref{lem:2-powers} shows that any element of order two in
the image of $\rho$ has vanishing trace.
Since any element of $\Gal(E/\Q)$ that is not of order two is a power of
$\sigma$,
we conclude from part~(ii) of the same lemma that
all the traces of $\rho$ lie in the ideal of $\F_2[X]/X^n$
generated by $\Trace(\rho(\sigma))$.
Since the trace of any element in the image of $\rhobar$ is zero,
we see that this ideal is contained in the maximal ideal of $\F_2[X]/X^n$.
Applying part~(iii) of Lemma~\ref{lem:2-powers}, we deduce
that $\Trace(\sigma)^{2^{m-1}} = 0$ (since $\sigma^{2^m} = 1$,
and so $\rho(\sigma^{2^m}) = I$).

Thus, on the one hand, the only
way that $X$ can arise as a trace of $\rho$ is if
$\Trace(\sigma) \equiv X \pmod X^2.$  
On the other hand, if this condition holds, then $X^{2^{m-1}} = 0,$
and hence $n \leq 2^{m-1}.$  This proves one direction of the ``if and only
if'' statement of the corollary.

Let us now prove the uniqueness assertion, assuming that
we are given a surjective map $R \rightarrow \F_2[X]/X^n$.
Since the corresponding triple $(V,L,\rho)$ 
deforms $(\Vbar,\Lbar,\rhobar),$ and since
$\sigma$ has non-trivial
image in $\Gal(\Q(\sqrt{-1})/\Q),$
while $\tau$ generates
the image of $I_N$ in $\Gal(E/\Q),$ 
we may choose a basis of
$V$ such that $\sigma$ and $\tau$ act through matrices in $\SL_2(\F_2[X]/X^n)$
of the form
$$
\rho(\sigma) = \left(\begin{matrix} a(\sigma) & b(\sigma) \\
c(\sigma) & d(\sigma) \end{matrix} \right) \equiv \left( \begin{matrix}
1 & 1 \\ 0 & 1\end{matrix} \right) \mod X \, ,
\quad
\rho(\tau) = \left(\begin{matrix} a(\tau) & 0 \\
c(\tau) & d(\tau) \end{matrix} \right) \equiv \left( \begin{matrix}
1 & 0 \\ 0 & 1\end{matrix} \right) \mod X.
$$
Now conjugating by matrices in
$\kernel(\GL_2(\F_2[X]/X^n) \rightarrow \GL_2(\F_2))$
of the form
$$\left(\begin{matrix} \alpha & 0 \\ \gamma & \delta \end{matrix}
\right),$$
it is easy to show that we may change our basis so that
$$\rho(\sigma) = \left(\begin{matrix} 1 + u X & 1 \\
u X & 1 \end{matrix} \right),
\qquad
\rho(\tau) = \left(\begin{matrix} 1 & 0 \\
u X & 1 \end{matrix} \right),
$$
for some $u \in (\F_2[X]/X^n)^{\times}.$  
Thus, after applying the inverse of the automorphism
of $\F_2[X]/X^n$ induced by
the map $X \mapsto u X$,
we see that we may put $\rho$ in the form
\begin{equation}\label{eq:normal form}
\rho(\sigma) = \left(\begin{array}{cc} 1 + X & 1  \\
X & 1 \end{array} \right),
\qquad
\rho(\tau) =  \left(\begin{array}{cc} 1  & 0  \\
X & 1 \end{array} \right).
\end{equation}
This proves the uniqueness statement.

Finally, one checks that the preceding formula
gives a well-defined homomorphism $\rho:\Gal(E/\Q)
\rightarrow \SL_2(\F_2[X]/X^n)$,
so long as $n\leq 2^{m-1}$,
and that it deforms $\rhobar$.  It is certainly flat at $2$,
since the inertia group at two acts through its image in
$\Gal(\Q(\sqrt{-1})/\Q)$.
Thus, if we let $L$ denote the line spanned by the vector $(0,1)$,
then we obtain an object of $\Def(\F_2[X]/X^n)$ of the required
sort (since $\Trace(\rho(\sigma)) = X$).  This completes
the proof of the corollary.
$\qed$
\end{Proof}

\

Let us consider the particular case $n=2$ of the preceding
corollary.

\begin{prop}\label{prop:tangent space at 2}
If $N\not\equiv 1 \pmod 8,$ then $\Def(\F_2[X]/X^2) = 0.$
If $N \equiv 1 \pmod 8,$ then $\Def(\F_2[X]/X^2)$ is one dimensional
over $\F_2.$  Furthermore, if $(V,L,\rho)$ corresponds to the
non-trivial element, then we have the following formulas for
the traces of $\rho$:
\begin{itemize}
\item[(i)] If $\ell$ is an odd prime distinct from $N$, then
$$\Trace(\rho(\Frob_{\ell}))
= \begin{cases}
0 \, \text{ if } \ell \equiv 1 \pmod 4 \text{ or } \ell \text{ is
a square} \mod N \\ X \text{ otherwise }  \end{cases}.$$
\item[(ii)]  If $\alpha_2$ denotes the eigenvalue of $\Frob_2$ on the
rank one $\F_2[X]/X^2$-module of $I_2$-coinvariants of $V$,
then
$$\alpha_2 = \begin{cases} 1 \text{ if } -4 \text{ is an $8$th power mod }
 N \\
1 + X \text{ if not }\end{cases}.$$
\end{itemize}
\end{prop}

\begin{Proof}
If $N \not\equiv 1 \pmod 8$ then Proposition~\ref{prop:2-cft}
shows that $H$ has order two, and
Corollary~\ref{cor:existence} shows that 
any map $R \rightarrow \F_2[X]/X^2$ factors through the 
map $R \rightarrow \F_2$.  Thus in this case $\Def(\F_2[X]/X^2) = 0,$
as claimed. 

If $N \equiv 1 \pmod 8,$ then conversely we conclude from
Proposition~\ref{prop:2-cft} that $H$ has order divisible by $4$.
Corollary~\ref{cor:existence} then shows that there is
a unique surjection $R \rightarrow \F_2[X]/X^2$, and
thus that $\Def(\F_2[X]/X^2)$ is one dimensional over $\F_2$.
If $F$ denotes the subextension of $E$ over $K$ cut out by
this non-trivial deformation, then $F$ is a dihedral extension
of $\Q$ of degree $8$, unramified over $K$, containing $K(\sqrt{-1})$.
(Concretely, as we saw in the proof
of Proposition~\ref{prop:2-cft}, the field $F$ has the form
$\Q(\sqrt{-1},\sqrt{\nu},\sqrt{\nubar}),$ for appropriate
$\nu,$ $\nubar$.)

Recall the presentation~(\ref{eq:presentation}) of $\Gal(E/\Q)$.
If we let $\sigmabar$ and $\taubar$ denote the image of $\sigma$
and $\tau$ under the surjection $\Gal(E/\Q) \rightarrow \Gal(F/\Q)$,
then $\Gal(F/\Q)$ has the following presentation:
$$\Gal(F/\Q)= \langle \sigmabar, \taubar | \sigmabar^4 = \taubar^2 =
(\sigmabar \taubar)^2 = 1 \rangle.$$
Recall from the proof of Corollary~\ref{cor:existence}
that the only elements of $\Gal(F/\Q)$ whose images under
$\rhobar$ have non-zero trace (which is then equal to $X$)
are $\sigmabar^{\pm 1};$ that is, the elements of $\Gal(F/\Q)$
that are of order $4$.  

If $\ell$ is an odd prime distinct from $N$, then 
$\ell$ is unramified in $F$.  The final remark of the preceding
paragraph shows that
$$\Trace(\rho(\Frob_{\ell})) = \begin{cases} 0 \text{ if } \Frob_{\ell}
\text{ has order $1$ or $2$}\\ X \text{ if } \Frob_{\ell} \text{ has
order $4$ } \end{cases}.$$
Now $K$ is the maximal subfield of $F$ fixed by $\sigmabar$,
while one checks that any element of $\Gal(F/\Q)$ of order two
fixes at least one of the subfields $\Q(\sqrt{-1})$ or $\Q(\sqrt{N})$
of $F$.  Thus we see that $\rho(\Frob_{\ell})$ has trace zero
(as opposed to trace $X$)
if and only if $\ell$ splits in at least one of the
fields $\Q(\sqrt{-1})$ or $\Q(\sqrt{N})$.  This establishes~(i).

Again referring to the presentation~(\ref{eq:presentation})
of $\Gal(E/\Q),$
one easily checks that $D_2(E/\Q)$ is generated
by $\sigma \tau$ and $\sigma^{2^{m-1}},$ with $I_2(E/\Q)$
being generated by $\sigma \tau$.  (Recall that 
$2^m$ denotes the order of $H$.)
Thus $D_2(F/\Q)$
is generated by $\sigmabar \taubar$ and $\sigmabar^{2^{m-1}}.$
Thus if $m \geq 3,$ then we see that $D_2(F/\Q) = I_2(F/\Q)$,
while if $m = 2$, then $D_2(F/\Q)/I_2(F/\Q)$ is generated by the
image of $\sigmabar^2$.
 
In terms of the explicit model~(\ref{eq:normal form})
for $\rho$, we see that the coinvariants of
$I_2(F/\Q) = \langle \sigmabar \taubar \rangle$ on $V$ are
spanned by the image of the basis vector $(0,1)$, and
that $\sigmabar^2$ (which is central in $\Gal(F/\Q)$,
and so does act on the space of coinvariants)
acts on the image of this vector as multiplication by
$1+X$.  Combining this computation with the discussion of the
previous paragraph proves part~(ii), once we recall
from Proposition~\ref{prop:2-cft} that $m\geq 3$ if and only
if $-4$ is an $8$th power modulo $N$.
$\qed$
\end{Proof}

\

The quotient $R/2$ is the universal deformation ring
classifying deformations of $(\Vbar,\Lbar,\rhobar)$
in characteristic $2$.  
The preceding two results together imply that
$R/2 \iso \F_2[X]/X^{2^{m-1}},$ where $2^m$ is
the order of $H$; formula~(\ref{eq:normal form})
then gives an explicit model for the universal deformation
over $R/2$.

\

We close this section by observing that 
Theorem~\ref{thm:main:p=2}
follows from Corollaries~\ref{EPR} and~\ref{cor:existence}
taken together.

\section{Explicit deformation theory: $p$ odd}\label{sec:odd explicit}

In this section we suppose that $p\geq 3$, and
that $N$ is prime to $p$.
We begin by considering the problem of analysing deformations
$(V,L,\rho) \in \Def(A)$, where~$A$
is an Artinian local $\F_p$-algebra with residue field $\F_p$.
Our results will be less definitive than those obtained
in the case of $p = 2$.

Let $\Delta$ denote the following
subgroup of $\GL_2(\F_p)\subset \GL_2(A)$:
$$\Delta =
\left \{
\left( \begin{matrix} \alpha & 0 \\ 0 & 1 \end{matrix}\right)
\, \left | \right.  \, \alpha \in \F_p^{\times} \right \},$$
and let $G'$ denote the kernel of the map
$\SL_2(A) \rightarrow \SL_2(\F_p)$
induced by reduction modulo $\m$ (the maximal ideal of $A$);
note that $G'$ is a
normal subgroup of $\GL_2(A)$.
If we let $G$ denote the subgroup of $\GL_2(A)$ generated by 
$G'$ and $\Delta$, then $G$ is isomorphic to the semi-direct product
$G' \sdp \Delta$, where $\Delta$ acts on $G'$ via conjugation.
Explicitly, one computes that
\begin{equation}\label{eq:conjugation}
\left( \begin{matrix} \alpha & 0 \\ 0 & 1 \end{matrix}\right)
\left( \begin{matrix} a & b \\ c & d \end{matrix} \right)
\left( \begin{matrix} \alpha^{-1} & 0 \\ 0 & 1 \end{matrix}\right)
= \left( \begin{matrix} a & \alpha b \\ \alpha^{-1} c & d \end{matrix} \right).
\end{equation}

\begin{lemma}\label{lemma:basis} Let $(V,L,\rho)$ be an object of $\Def(A)$,
and let $M$ denotes the finite flat group scheme over $\Z_p$
whose generic fibre equals $V$.  If
$0 \rightarrow M^{\co} \rightarrow M \rightarrow M^{\et} \rightarrow 0$
denotes the connected-\'etale exact sequence of $M$, 
then there is a basis for $V$ over $A$ such that

\begin{itemize}
\item[(i)] The representation $\rho: G_{\Q} \rightarrow \GL_2(A)$
has image lying in $G$.

\item[(ii)]   The submodule $M^{\co}(\Qbar_p)$ of $V$
(which is free of rank one, by Corollary~\ref{uniquenesscor})
is spanned by the vector $(1,0)$.

\item[(iii)] The line $L$ is spanned by the vector $(1,1)$.

\end{itemize}
\end{lemma}

\begin{Proof} By definition of the deformation problem
$\Def,$ the determinant of $\rho$ is equal to $\overline{\chi}_p$.
Thus $\image(\rho)$ sits in the exact sequence of groups
$$0 \rightarrow G' \rightarrow \image(\rho) \rightarrow \F_p^{\times}
\rightarrow 0.$$
The order of $\F_p^{\times}$ is coprime to the order of $G'$,
and so this exact sequence splits.
If we fix a splitting $s$,
then one easily sees that we may choose an eigenbasis for
the action of $s(\F_p^{\times})$ so that this group acts via
the matrices in $\Delta$.   Thus condition~(i) is satisfied
for this basis.  Condition~(ii) follow directly from condition~(i).
The stipulations of the deformation problem $\Def$ then imply that
$L$ is spanned by a vector of the form $(1,u)$, for some unit
$u \in A^{\times}$.  Rescaling the second basis vector by $u$,
we may assume that $L$ is in fact spanned by $(1,1)$.
$\qed$
\end{Proof}

\

From now on, we fix an object $(V,L,\rho) \in \Def(A)$,
and choose a basis of $V$ as in the preceding lemma.
Thus we may regard $\rho$ as a homomorphism
$G_{\Q} \rightarrow G \subset \GL_2(A).$

\begin{lemma}\label{lem:inertia at N}
If $(V,L,\rho)$ is a non-trivial deformation,
then the image of 
$I_N$ under $\rho$ is a cyclic subgroup of $G'$ of order $p$.
Furthermore,
if $\left(\begin{matrix} a & b \\ c & d \end{matrix}\right)$
is a generator of this cyclic group, then neither $b$ nor $c$
vanishes, and neither $a$ nor $d$ equals $1$.
\end{lemma}

\begin{Proof}
Since the image under $\rho$
of inertia at $N$ acts trivially on each of the lines $L$ and $V/L$
(the determinant of $\rho$ equals $\overline{\chi}_p,$
which is trivial on $I_N$),
we see that $I_N$ acts via an abelian group of exponent~$p$. Since
tame inertia is pro-cyclic, inertia at $N$ must act
through a group of order dividing~$p$. 
If $I_N$ has trivial image, then Proposition~\ref{mindefiso}
shows that $V = A\otimes_{\F_p} \Vbar,$ and thus
that $(V,L,\rho)$ is the trivial deformation, contradicting
our assumption.  Thus $I_N$ has image of order $p$.

The line $L$ is spanned by the vector $(1,1)$.
Thus 
if $\gamma = \left(\begin{matrix} a & b \\ c & d \end{matrix}\right)$
is a generator of the image of $I_N$, it fixes such a vector.
Since the determinant of $\rho$ equals $\overline{\chi}_p,$ we see
that $\det(\gamma) = 1$.
If $a$ (respectively $d$) equals $1$ then 
we conclude that $b$ (respectively $c$) equals $0$.
If either $b$ or $c$ vanishes, one easily checks that $\gamma$
must be the identity, contradicting the fact that $I_N$ has
non-trivial image. 
$\qed$
\end{Proof}

\

If $F$ denotes the extension of $\Q$ cut out by the
kernel of $\rho,$  then $F$ contains $\Q(\zeta_p)$ (where
$\zeta_p$ denotes a primitive $p$th root of unity),
since $\det(\rho) = \overline{\chi}_p$.  
We let $F^{ab}$ denote the maximal subextension
of $F$ abelian over $\Q(\zeta_p)$.

\begin{lemma} \label{lemma:classf} 
The $p$-part of the 
conductor of $F^{ab}/\Q(\zeta_p)$ divides $\pi^2$, where
$\pi  = (1-\zeta_p) \Z[\zeta_p]$, and the extension $F/\Q(\zeta_p)$ has
inertial degree dividing $p$  
at $N$ and is unramified outside $N$ and $\pi$. \end{lemma}

\begin{Proof}
Lemma~\ref{lem:inertia at N} shows that
the image under $\rho$ of inertia at $N$ is a cyclic
group of order dividing $p$.
Therefore it suffices to prove the conductor bound at $\pi$.

The image under $\rho$
of $G_{\Q(\zeta_p)}$ lies in $G'$, a $p$-group, and so we see
that $\Gal(F^{ab}/\Q(\zeta_p))$ is an abelian $p$-group.
Thus it is a compositum of cyclic
extensions of $p$-power degree.  The conductor of a compositum
of cyclic extensions is equal to the g.c.d.~of the 
conductors of the individual cyclic extensions, and
thus it suffices to bound the conductor of a cyclic
subextension of $F^{ab}$ of degree $p^k$, for some $k\geq 1$.

Let $F'$ be such a subextension, and suppose that
the conductor of $F'$ is divisible by $\pi^3$.
There are $(p-1) p^{k-1}$ faithful characters of $\Z/p^k$,
and so by the conductor discriminant formula,
the discriminant $\Delta_{F'/\Q(\zeta_p)}$ is divisible by
$\pi^{3(p-1)p^{k-1}}$.  Thus
the $p$-root
discriminant of $F'$ satisfies
$$\delta_{F',p}  \ge \delta_{\Q(\zeta_p)} 
N_{\Q(\zeta_p)/\Q}(\pi^{3(p-1)p^{k-1}})^{1/[F':\Q]}
= p^{(p-2)/(p-1)} \cdot p^{3(p-1)/(p(p-1))}$$
and thus
$$v_p(\D_{F'/\Q}) \ge 1 + \frac{1}{p-1} + \frac{p-3}{p(p-1)}.$$
This violates Fontaine's bound \cite{fbound} when $p \ge 3$.
The result follows for $F^{ab}$.
$\qed$
\end{Proof}

\

In order to apply this result, we will need to classify the
relevant class fields of $\Q(\zeta_p)$ that can arise in
the situation of the preceding lemma.

\begin{prop}\label{prop:odd cft}
Let $p$ be an odd prime, and let $N$ be a prime distinct from $p$.
For any value of $i$, let $K_{(i)}$ denote the maximal abelian extension
of $\Q(\zeta_p)$ satisfying the following conditions: $K_{(i)}$ has conductor
dividing $\pi^2 N$; the Galois group $\Gal(K_{(i)}/\Q_p(\zeta_p))$ has exponent
$p$;
the Galois group $\Gal(\Q_p(\zeta_p)/\Q)$ acts on $\Gal(K_{(i)}/\Q_p(\zeta_p))$
through the $i$th power of the mod $p$ cyclotomic character $\chibar_p$.
Then:


\begin{itemize}

\item[(i)] $K_{(1)} = \Q(\zeta_p,N^{1/p})$;

\item[(ii)] $K_{(0)}  = \begin{cases} \text{ the degree $p$ subextension
of $\Q_p(\zeta_p, \zeta_N)/\Q_p(\zeta_p)$ if } N \equiv 1 \mod p \\
\Q(\zeta_p) \text{ otherwise }
\end{cases};$

\item[(iii)]
$K_{(-1)} = \begin{cases} \text{ a degree $p$ extension of
$\Q_p(\zeta_p)$ if } N^2 \equiv 1 \mod p \\ \Q(\zeta_p)
\text{ otherwise }
\end{cases}.$
\end{itemize}
\end{prop}

\begin{Proof}
Let $E_{(i)}$ denote the unramified extension of
$\Q(\zeta_p)$ of exponent $p$ corresponding to
the maximal elementary $p$-abelian quotient of
the class group of $\Q(\zeta_p)$ on which
$\Gal(\Q(\zeta_p)/\Q)$ acts through $\chibar_p^i$.
Then we have the short exact sequence of abelian
Galois groups
$$
0 \rightarrow \Gal(K_{(i)}/E_{(i)}) \rightarrow
\Gal(K_{(i)}/\Q(\zeta_p)) \rightarrow \Gal(E_{(i)}/\Q(\zeta_p))
\rightarrow 0.$$
Global class field theory allows us to compute the group
$\Gal(K_i/E_{(i)})$.
Indeed, it sits in the exact sequence
$$ (\Z[\zeta_p])^{\times}
\longrightarrow \left( \left(
\Z[\zeta_p]/\pi^2 \times \Z[\zeta_p]/N \right)^{\times} / p \right)_{(i)}
\longrightarrow \Gal(K_{(i)}/E_{(i)}) \longrightarrow 0;$$
here the subscript $(i)$ denotes the maximal quotient
on which $\Gal(\Q(\zeta_p)/\Q)$ acts via $\chibar_p^i$.

Since the reduction mod $\pi^2$ of 
the global unit $\zeta_p = 1 + (\zeta_p -1)$ generates the $p$-power
part of $(\Z[\zeta_p]/\pi^2)^{\times}$, we may eliminate
this factor from the second term of the preceding exact sequence.
If we fix a prime $\frak n$ over $N$ in $\Z[\zeta_p]$,
then as in the proof of Lemma~\ref{four},
we obtain a surjection
$$(\Z[\zeta_p]/\frak n)^{\times}/(p,N^{1-i} - 1) \rightarrow
\Gal(K_{(i)}/E_{(i)}).$$
Consequently, we find that
$\Gal(K_{(i)}/E_{(i)})$ is either trivial
(when $N^{(1-i)} \not\equiv 1 \pmod p$) or of order $p$
(when $N^{(1-i)} \equiv 1 \pmod p$).

Let us now consider the particular cases $i  = 1,0,-1.$
The $1$, $0$ and $-1$ eigenspaces inside the class group of
$\Q(\zeta_p)$  are trivial by Kummer theory, 
abelian class field theory and Herbrand's theorem respectively.
Thus for these values of $i$, we have $E_i = \Q(\zeta_p)$,
and so the preceding paragraph yields a computation of
$\Gal(K_{(i)}/\Q(\zeta_p))$.  The explicit descriptions
of $K_{(i)}$ in the case when $i = 1$ or $0$ are easily
verified, and so we leave this verification to the reader.
$\qed$
\end{Proof}

\

We are now in a position to determine the reduced Zariski tangent space
to the deformation functor $\Def$.  We will also record
some useful information regarding non-trivial elements of this tangent
space (assuming that they exist).

\begin{prop}\label{prop:tangent space at odd p}
If $p$ does not divide the numerator of $(N-1)/12$,
then $\Def(\F_p[X]/X^2) = 0$; otherwise,
$\Def(\F_p[X]/X^2)$ is one dimensional over $\F_p.$
Suppose for the remainder of the statement of the proposition
that we are in the second case,
and let 
$(V,L,\rho)$ correspond to a
non-trivial element of $\Def(\F_p[X]/X^2)$.

(i) If as above $F$ denotes the extension of $\Q$ cut out 
by the kernel of $\rho$, then
$F$ is equal to the compositum $K_{(1)}K_{(0)}K_{(-1)}$
$($where the class fields $K_{(i)}$ of $\Q(\zeta_p)$
are defined as in the statement of the previous 
proposition$)$.

(ii) If $p = 3$, then $\Gal(F/\Q(\zeta_p)) \iso
\Gal(K_1/\Q(\zeta_p)) \times \Gal(K_0 /\Q(\zeta_p)),$
and the image of an appropriately chosen generator of the first $($respectively
second$)$ factor under $\rho$ has the form
$\left(\begin{matrix} 1 & -r X \\ r X & 1\end{matrix} \right)$
$($respectively $\left(\begin{matrix} 1 + r X &  0 \\  0 & 1 - r X \end{matrix} \right))$
for some $r \in \F_p^{\times}$.

(iii) If $p \geq 5$, then $\Gal(F/\Q(\zeta_p)) \iso
\Gal(K_1/\Q(\zeta_p)) \times \Gal(K_0 /\Q(\zeta_p)) \times \Gal(K_1/\Q(\zeta_p)),$
and the image of an appropriately chosen generator of the first $($respectively
second, respectively third$)$ factor under $\rho$ has the form
$\left(\begin{matrix} 1 & -r X \\ 0 & 1\end{matrix} \right)$
$($respectively $\left(\begin{matrix} 1 + r X &  0 \\  0 & 1 - r X \end{matrix} \right),$
respectively
$\left(\begin{matrix} 1 &  0 \\ r X & 1\end{matrix} \right))$
for some $r \in \F_p^{\times}$.

(iv)
We have the following formulas for
the traces of $\rho$:
\begin{itemize}
\item[(iv.i)] If $\ell$ is a prime distinct from $N$ and $p$,
then $$\Trace(\rho(\Frob_{\ell}))
= \begin{cases}
1 + \ell \, \text{ if } \ell \equiv 1  \pmod p \text{ or } \ell \text{ is
a $p$th power} \mod N \\  1 + \ell + u X \text{ otherwise }  \end{cases};$$
here $u$ denotes an element of $\F_p^{\times}$.

\item[(iv.ii)]  If $\alpha_p$ denotes the eigenvalue of $\Frob_p$ on the
rank one $\F_p[X]/X^2$-module of $I_p$-coinvariants of $V$,
then
$$\alpha_p = \begin{cases} 1 \text{ if } p \text{ is a $p$th power mod }
 N \\
1 + u X \text{ if not }\end{cases};$$
again, $u$ denotes an element of $\F_p^{\times}$.
\end{itemize}
\end{prop}

\begin{Proof}
Let $(V,L,\rho)$ be a non-trivial element of $\Def(\F_p[X]/X^2)$,
cutting out the extension $F$ of $\Q$.  As above, we choose the basis
of $V$ so that the conditions of Lemma~\ref{lemma:basis} are satisfied.
Since $G'$ is abelian,
we see that $F = F^{ab}$.  Equation~\ref{eq:conjugation}, together
with Lemma~\ref{lemma:classf}, thus shows
that $F$ is contained in the compositum $K_{(1)}K_{(0)}K_{(-1)}.$
Lemma~\ref{lem:inertia at N} then shows that in fact
$F$ must be equal to this compositum, proving part~(i) of the proposition;
that furthermore,
each of the extensions $K_{(1)},$ $K_{(0)}$ and $K_{(-1)}$ of
$\Q(\zeta_p)$ must be non-trivial, and thus that
$N \equiv 1 \mod p$, by proposition~\ref{prop:odd cft};
and that either part~(ii) or part~(iii)
of the proposition is satisfied, depending 
on whether $p = 3$ or $p \geq 5$.  (We choose the generator of
each group $\Gal(K_{(i)}/\Q(\zeta_p))$ to be the image of some
fixed generator of the inertia group $I_N$.)

Suppose conversely that 
$N \equiv 1 \mod p$, so that each of $K_{(1)},$ $K_{(0)}$ and
$K_{(-1)}$ is a non-trivial extension of $\Q(\zeta_p)$.
Write $F = K_{(1)}K_{(0)} K_{(-1)}$.
If we fix an element $r \in \F_p^{\times},$ then
we may use the formulas of parts~(ii) and~(iii) to define 
a representation $\rho:\Gal(F/\Q) \rightarrow G \subset \GL_2(\F_p[X]/X^2)$.
If we let $L$ denote the line spanned by $(1,1)$, then this representation
will deform the representation $(V,L)$.
Thus it will provide an element of $\Def(\F_p[X]/X^2)$ provided that
it is finite at $p$.  An argument as in the proof of Lemmas~\ref{three}
and~\ref{four} shows that this is automatically the case when
$p \geq 5,$ and holds provided $p$ divides $(N-1)/12$, when $p = 3$.
This establishes the initial claim of the proposition.

It remains to prove part~(iv) of the proposition.
Suppose first that $\ell$ is a prime distinct from $p$ and $N$.
We may write 
$$\rho(\Frob_{\ell}) = \left( \begin{matrix} \ell & 0 \\ 0 & 1 \end{matrix}
\right) \left( \begin{matrix} 1 + a X & b X \\ c X & 1 - a X \end{matrix}
\right),$$
for some elements $a,b,c \in \F_p$.  Thus
$\Trace(\rho(\Frob_{\ell})) = 1 + \ell + (\ell -1 ) a X.$
This is distinct from $1 + \ell$ if and only if $\ell \not\equiv 1 
\mod p$ and $a \neq 0$.  The latter occurs if and only
if the primes over $\ell$ are not split in the extension
$K_{(0)}/\Q(\zeta_p),$ which in turn is the case if and only if
$\ell$ is not a $p$th power $\mod N$.  (Here we have taken into
account the explicit description of $K_{(0)}$ provided by
Proposition~\ref{prop:odd cft}.)  Thus we proved part~(iv.i).

Since the vector $(1,0)$ spans the subspace
$M^{\co}(\Qbar_p)$ of $V$, the space of $I_p$-inertial coinvariants
is spanned by the image of the vector $(0,1)$.  The Frobenius element
$\Frob_p$ acts non-trivially on the image of this vector if and only
if the prime over $p$ is not split in the extension $K_{(0)}$ of
$\Q(\zeta_p)$, which is the case if and only if $p$ is not a $p$th
power $\mod N$.  This proves part~(iv.ii).
$\qed$
\end{Proof}

\

As we will see below,
for $p \ge 3$, the rank $g_p + 1$ of $\T/p$ over $\F_p$  is no
longer explained by an abelian extension of number fields
(and hence by a single class group), as it is in the case $p = 2$,
but by certain more complicated solvable extensions.
However, the question of whether or not
$g_p=1$ is somewhat tractable.
Indeed, from Corollary~\ref{EPR} we deduce the following criterion.

\begin{lemma}\label{lem:criterion}
The rank $g_p$  of the parabolic Hecke algebra
$\T^0/p$ over $\F_p$ is greater than one $($equivalently, $\displaystyle
\T^0 \neq Z_p)$
if and only if there
exists a $(V,L,\rho)$ in $\Def(\F_p[X]/X^3)$ whose
traces generate $\F_p[X]/X^3$. \end{lemma}

In order to apply this lemma, we now assume that $A = \F_p[X]/X^3,$
so that $(V,L,\rho)$ lies in $\Def(\F_p[X]/X^3)$.  As always,
we assume that the basis of $V$ is chosen so as to satisfy the
conditions of Lemma~\ref{lemma:basis}.
We let 
$\rho_n$ denote the composition of $\rho$ with
the natural surjection 
$\GL_2(\F_p[X]/X^3) \rightarrow \GL_2(\F_p[X]/X^n)$,
for $n \leq 3$.
Requiring the traces of $\rho$ to generate $\F_p[X]/X^3$
is equivalent to requiring the traces of $\rho_2$ to
generate $\F_p[X]/X^2,$ which in turn is equivalent
is to requiring that $\rho_2$ be a non-trivial deformation.
We assume this to be the case.
Also, we let $F_n$ denote the extension cut out
by the kernel of $\rho_n$. Thus
$F_1 =\Q(\zeta_p)$, and $F_3 = F$. 

Since we are assuming that $\rho_2$ is non-trivial,
Proposition~\ref{prop:tangent space at odd p}
shows that $p$ divides the numerator of $(N-1)/12$,
and that $F_2$ is equal to the compositum of the
class fields $K_{(i)}$ (for $i = 1,0,-1$).

\begin{lemma}\label{lemma:shape of F_2}
(i) We have $F_2 = F^{ab}$, and $\Gal(F_2/F_1) \iso (\Z/p)^2$
$($respectively $(\Z/p)^3)$ if $p = 3$ $($respectively $p \geq 5)$.

(ii) $F/F_2$ is unramified at $N$.
\end{lemma}

\begin{Proof}
Since $p \geq 3,$ we see that $G'$ has exponent $p$.
Lemma~\ref{lemma:classf} and equation~(\ref{eq:conjugation})
then imply that
$F^{ab} \subset K_{(1)}K_{(0)} K_{(-1)} = F_2$.
Certainly $F_2 \subset F^{ab}$, and so we have the equality stated
in~(i).
The claims regarding $\Gal(F_2/F_1)$ follow from
parts~(ii) and~(iii) of Proposition~\ref{prop:tangent space at odd p}.

Part~(ii) follows from the Lemma~\ref{lem:inertia at N}
and the fact that $F_2/F_1$ is ramified at $N$.
$\qed$
\end{Proof}

\

We now separate our analysis into two cases:  $p = 3$, and $p \geq 5$.

%

\subsection{$p = 3$}

Throughout this subsection we set $p = 3$.

\begin{lemma} The extension $F/F_2$ is unramified everywhere
and has degree exactly three.
\end{lemma}

\begin{Proof}
The image of
$\rho_{| G_{\Q(\sqrt{-3})}}$
is a subgroup of
$G' = \ker(\GL_2(\F_3[X]/X^3) \rightarrow \GL_2(\F_3))$
whose image in $\GL_2(\F_3[X]/X^2)$ is isomorphic to
$(\Z/3)^2$, by Lemma~\ref{lemma:shape of F_2}.  Thus the commutator
subgroup of the image of 
$\rho_{| G_{\Q(\sqrt{-3})}}$
is either trivial or cyclic of order three.
Thus the extension $F/F_2$ has degree at most three.

Consider the representation $\rho_2$, which factors through
$\Gal(F_2/\Q)$.  By assumption this yields a non-trivial element of
$\Def(\F_3[X]/X^2)$.
Part~(ii) of Lemma~\ref{prop:tangent space at odd p} thus shows
that the image under $\rho_2$ of the element of order three coming
from the $\chibar_p^1$ extension $K_{(1)} = \Q(\sqrt{-3},\sqrt[3]{N})$
of $\Q(\sqrt{-3})$
must be of the form
$$\left(\begin{matrix} 1 & 
 -r X \\  r X & 1 \end{matrix} \right),$$
and that the image under $\rho_2$ of
the element of order three coming from the
$\chibar_p^{0}$ extension $K_{(0)}$ 
of $\Q(\sqrt{-3})$ is of the form
$$\left(\begin{matrix} 1+rX & 0 \\ 0 & 1-rX \end{matrix} \right),$$
for some $r \in \F_3^{\times}$.
Lifting these two elements (in any way) to $\GL_2(\F_3[X]/X^3)$ and
taking their commutator, we produce a new element in $\Gal(F/\Q)$
which has a lower left-hand entry equal to $r^2 X^2 = X^2$.
This element cannot be in the decomposition group at $3$ because
it doesn't preserve  $M^0(\Qbar_3)$, which is generated by
$(1, 0)$. Thus
$F/F_2$ has order exactly three and is unramified at
all primes above three.  Part~(ii) of Lemma~\ref{lemma:shape of F_2}
shows that the extension $F/F_2$ is also unramified at all
primes above $N$, and the lemma is proved.
$\qed$
\end{Proof}

\

Let $K = \Q(\sqrt[3]{N})$, and as above write $K_{(1)}
= K^{gal} = K(\sqrt{-3})$.
The extension $F/K_{(1)}$ has degree $9$, and $\Gal(F/K_{(1)}) = (\Z/3\Z)^2$.
Moreover, $F/K_{(1)}$ is unramified everywhere. The following lemma
shows that the existence of such an extension $F$
is sufficient for the construction of a deformation $\rho$
of the type considered here.
This completes the proof of part one of each of
Theorems~\ref{theorem:oddp} and~\ref{theorem:merel}.

\begin{lemma} \label{lemma:hwl}
If $N \equiv 1 \mod 9$, then
the class group  of $K_{(1)}=\Q(\sqrt{-3},\sqrt[3]{N})$ has
$3$-rank greater than or equal \emph{(}equivalently, equal\emph{)} to 
two if and only if there exists a surjection
$R \rightarrow \F_3[X]/X^3$; the kernel of the corresponding deformation
$\rho:G_{\Q} \rightarrow \GL_2(\F_p[X]/X^3)$ then
cuts out the $(3,3)$ unramified class field $F$ of $K_{(1)}$.
\end{lemma}

\begin{Proof}
The preceding discussion establishes the ``if'' claim, and so
it suffices to prove the ``only if'' claim. 
Genus theory and a consideration of
the ambiguous class predicts that the 3-rank of the class group
of $K_{(1)}$ is either
one or two, and hence by assumption this rank
is exactly two (see for example \cite{Gerth}).
We let $F$ denote the corresponding
unramified $(3,3)$-extension of $K_{(1)}$,
and (as above) let $F^{ab}$ denote the
unique subextension of $F$ abelian over $\Q(\sqrt{-3})$.
It is easily checked that $F^{ab}$ is in fact
the maximal abelian 3-power extension of $\Q(\sqrt{-3})$
that is unramified over $K_{(1)}$,
and that $F^{ab} = K_{(1)}K_{(0)}$.

Proposition~\ref{prop:tangent space at odd p} yields a Galois
representation $\Gal(F^{ab}/F) \rightarrow \GL_2(\F_3[X]/X^2)$,
while lemma~\ref{lem:cf structure} below shows
that $\Gal(F/\Q(\sqrt{-3}))$ is the unique non-abelian group of
order $27$ and of exponent three.
It is then easy to see that one can lift the representation of $\Gal(F^{ab}/F)$
to a representation $\rho: \Gal(F/\Q) \rightarrow \GL_2(\F_3[X]/X^3).$
Furthermore one checks that for any such lift, the image of $I_N$
fixes an appropriate line.

To show that we have constructed an element of $\Def(\F_3[X]/X^3)$,
as required, it remains to 
show that this representation 
extends to a finite flat group scheme at $3$,
For this, it suffices to work over the maximal unramified extension
of $\Q_3$. Since $F/\Q(\sqrt{-3})$ is unramified at $3$
(because $N \equiv 1 \mod 9$),
the representation $\rho |_{\Q^{ur}_3}$ factors through
a group of order two, and explicitly prolongs to a product of
trivial and multiplicative group schemes.  Thus $\rho$ is indeed
finite at the prime $3$.
$\qed$
\end{Proof}

\begin{lemma}\label{lem:cf structure}
The Galois group $\Gal(F/\Q(\sqrt{-3}))$ is the (unique up to isomorphism)
non-abelian group of
order $27$ and of exponent three.
\end{lemma}

\begin{Proof}
Let $\Gamma = \Gal(K_{(1)}/\Q(\sqrt{-3})) = \langle \gamma \rangle$.
The $3$-class group $H$ of $K_{(1)}$ is naturally a $\Z_3[\Gamma]$-module.
From class field theory we have that $H/(\gamma-1)H$ is isomorphic to
the Galois group over $K_{(1)}$
of the maximal abelian 3-extension of $\Q(\sqrt{-3})$
that is unramified over $K_{(1)}$; that is, to $\Gal(F^{ab}/K_{(1)})$,
a cyclic group of order 3. 
Thus by Nakayama's lemma $H$ is a cyclic $\Z_3[\Gamma]$-module.
By class field theory, the quotient $H/3$ is isomorphic to $\Gal(F/K_{(1)})$.

Note that $\Gal(F/\Q(\sqrt{-3}))$ sits in the exact sequence:
$$0 \rightarrow \Gal(F/K_{(1)}) \longrightarrow
\Gal(F/\Q(\sqrt{-3}) \longrightarrow \Gal(K_{(1)}/
\Q(\sqrt{-3})) \longrightarrow 0,$$
which is an extension of $\Gamma \cong \Z/3\Z$ by $H/3 \cong (\Z/3\Z)^2$.
The action via conjugation of $\Gamma$
on $H/3$ is non-trivial, since otherwise $H$ could
not be cyclic as a $\Gamma$-module. Already this shows that
$\Gal(F/\Q(\sqrt{-3}))$ is one of the two non-abelian groups of
order $27$.
To pin down the group precisely, we must show that it has
exponent three.  For this, it suffices to find a
splitting of the above exact sequence (a section from $\Gamma = \Gal(K_{(1)}/
\Q(\sqrt{-3}))$ back to $\Gal(F/\Q(\sqrt{-3}))$).
Since the inertia group above $N$ in $\Gal(F/\Q(\sqrt{-3}))$ has
order exactly three, and maps isomorphically to 
$\Gamma$, the required splitting exists.
$\qed$
\end{Proof}

\

The final result of this section provides a relation between
the rank of the $3$-class group of $K_{(1)}$ and the power of $3$
dividing the class number of $K$.

\begin{lemma} \label{lemma:ump}
The $3$-class group of $K_{(1)} = \Q(\sqrt{-3},\sqrt[3]{N})$ has
three rank two if the $3$-class group of $K = \Q(\sqrt[3]{N})$  $($which
is cyclic$)$ is divisible by nine.
\end{lemma}

\begin{Proof}
One has a class number relation between
$K$ and $K_{(1)}$ given by
$h_{K_{(1)}} = h^2_K/3 \cdot q$, where $q$ is the index of
the units in $K_{(1)}$ coming from $K$, $K^{\gamma}$, and $\Q(\sqrt{-3})$
inside the full unit group. (Here, as above, $\gamma$ denotes a generator
of the cyclic group $\Gamma = \Gal(K_{(1)}/\Q(\sqrt{-3}))$.) 
If $9 | h_K$, then
$27 | h_{K_{(1)}}$.  Recall from the proof of the previous lemma that 
the $3$-part $H$ of the class group is a cyclic $\Z_3[\Gamma]$-module,
and satisfies the condition that $H/(\gamma - 1)H$ is cyclic of order 3.

Now $\Z_3[\Gamma]$ 
admits no quotients $H'$ that are cyclic groups of order $27$ with the property
that $H'/(\gamma - 1)H'$ is of order $3$.
It follows that if $H$ is of order divisible by $27$, then it must be
non-cyclic, as claimed.
$\qed$
\end{Proof}

\

We conjecture that the converse to the preceding lemma is also true.
To prove this, it would suffice to show that whenever
$3 \| h_K$, the unit index $q$ is always equal to one. We have
verified this for all primes less than $50$,$000$ for which $3 \| h_K$.

\subsection{$p \ge 5$} 

Throughout this section we assume that $p \geq 5,$ and that
we are given a deformation to $\F_p[X]/X^3$ as in the 
discussion following Lemma~\ref{lem:criterion}.
Proposition~\ref{prop:tangent space at odd p} and Lemma~\ref{lemma:shape of F_2}
together show that $F_2 = K_{(1)} K_{(0)} K_{(-1)},$
that $\Gal(F_2/F_1) = (\Z/p)^3$, and that $F_2 = F^{ab}$.
We see that $F_2/F_1$ is unramified at $p$ if and
only if $N \equiv 1 \mod p^2$.

It follows from our determination of $F_2$ that $\Gal(F/F_1)$
is the full kernel of the map from 
$\SL_2(\F_p[x]/x^3)$ to $\SL_2(\F_p)$, since all the elements
of  
$$\Ker(\SL_2(\F_p[x]/x^3) \rightarrow \SL_2(\F_p[x]/x^2))$$
are generated by commutators of lifts of elements of
 $\Ker(\SL_2(\F_p[x]/x^2) \rightarrow \SL_2(\F_p))$.

\

\begin{lemma} \label{lemma:messy} If $E$ is a degree $p$ 
Galois extension of $F_2$ inside $F_3$ on which the matrix
$$\left(\begin{matrix} 1 + x^2 & 0 \\ 0 & 1 - x^2 \end{matrix}\right)$$
acts non-trivially, then $E/F_2$
is everywhere unramified.
\end{lemma}

\begin{Proof} Part~(ii) of Lemma~\ref{lemma:shape of F_2} shows
that this extension is unramified at primes above $N$.
To see that it is unramified at primes above $p$, it
suffices to note that the matrix 
$\left(\begin{matrix} 1 + x^2 & 0 \\ 0 & 1 - x^2 \end{matrix}\right)$
does not fix the vector $(1,0)$ (which spans $M^0(\Qbar_p)$).
$\qed$ \end{Proof}

\

Let $K = \Q(N^{1/p})$ and $L = K^{gal} = K(\zeta_p)$.

\begin{lemma}\label{lemma:final} The Hilbert class field of $K$
has $p$-rank at least two.
\end{lemma}

\begin{Proof}
Let us  first consider the extension $\Gal(F/K)$.
One sees that
$\Gal(F/K)^{ab} \iso (\Z/p \Z)^2\times (\Z/p)^{\times}$ is explicitly generated by the images
of
$$\left(\begin{matrix} 1 + x^k & 0 \\ 0 & (1+x^k)^{-1} \end{matrix}\right),$$ 
for $k=1,2$, together with the image of $\Delta$.
We let $H$ be the $(p,p)$-extension of $K$ contained in $F$ that is fixed by
$\Delta$.  We will show that $H$ is unramified over $K$.

We may write $H$ as a compositum $H = H_1 H_2$, where 
for each of $k=1,2$, we let $H_k$ denote a $p$-extension of $K$ contained in $F$,
on which the matrix
$\left(\begin{matrix} 1 + x^k & 0 \\ 0 & (1+x^k)^{-1} \end{matrix}\right)$ 
acts non-trivially.
If we let $\zeta^+_N$ denote an element of $\Q(\zeta_N)$ that
generates the degree $p$ subextension over $\Q$ (so that
$K_{(0)} = \Q(\zeta_p,\zeta_N^+)$),
then we may take $H_1$ to be
$K(\zeta^+_N)$, which is clearly
unramified everywhere over $K$ (it is the genus field).
We will show that $H_2$ is is also unramified
everywhere over $K$.
Lemma~\ref{lem:inertia at N} takes care of the primes above $N$,
and so it remains to treat the primes above $p$.

We begin by proving that $H_2(\zeta_p)/L$ is unramified. 
Lemma~\ref{lemma:messy} shows that
the extension $H_2 \cdot F_2/F_2$ is unramified.
Since $F_2/L$ is unramified, it follows that $H_2(\zeta_p)/L$
is unramified, as claimed.
We now use the fact that $H_2(\zeta_p)/L$
is unramified to show that $H_2/K$ is unramified. 
We consider two cases.
Suppose first that $p \| N - 1.$
Then $K$ is totally ramified at $p$,
and thus if $H_2/K$ is ramified
we deduce that since $H_2$ is \emph{Galois} over $K$,
$e_p(H_2) = p^2$, contradicting the fact that
$H_2(\zeta_p)/L$ is unramified.
If instead $N \equiv 1 
\mod p^2$, then things are even easier: If
$H_2/K$
is ramified at at least one prime $\frak{p}$ above $p$,
then again using the fact that $H_2/K$ is
Galois we deduce that
$p | e_{\frak{p}}(H_2)$. Yet $p$ is tamely ramified
in $L$ and therefore also in $H_2(\zeta_p)$. Thus
$H_2/K$ is unramified everywhere, and
$K$ has $p$-rank at least two. $\qed$
\end{Proof}

\

This completes 
the proof of parts two of Theorem~\ref{theorem:oddp} 
and Theorem~\ref{theorem:merel}.
We expect (based on the numerical evidence) that the condition that
the class group of $K$ has $p$-rank two
is equivalent
to the existence of an appropriate group scheme,
and thus to $g_p > 1$. 
Part of this could perhaps be proved by
more sophisticated versions of Lemmas~\ref{lemma:hwl},
\ref{lem:cf structure},
and~\ref{lemma:ump}.

\section{Examples}
The first example in Mazur's
table \cite{eisenstein} where $e_2 > 1$ occurs when
$N = 41$. The class group of $\Q(\sqrt{-41})$ is
$\Z/8\Z$. Thus one has $e_2 = 3$.
Using {\tt gp} and William Stein's programmes one can verify
that the class group of  $\Q(\sqrt{-21929})$  is
$\Z/256\Z$ and that $e_2 = 127$ for $N = 21929$.
In Mazur's table, $e_3$ always equals $1$ or $2$.
One has to go quite some distance before finding
an example where $e_3 > 2$.
For $N = 2143$, however, one has $e_3 = 3$. This 
is related to the fact that $2143$ is the smallest prime
congruent
to $1 \mod 9$
such that the class group of the corresponding extension
$K_{(0)}$ of $\Q(\sqrt{-3})$ (in the terminology of
Proposition~\ref{prop:odd cft})
has an element of order $9$. The corresponding class field
contributes to the maximal unramified \emph{solvable} extension
of $K = \Q(\sqrt[3]{2143})$. Finally, let us note
that when $p=3$, Lemmas~\ref{lemma:hwl} and~\ref{lemma:ump}
show that the value of $g_p$ is related to the
\emph{size} of the $3$-power part of the class group of $\Q(\sqrt[3]{N})$,
whereas for $p \ge 5$, Lemma~\ref{lemma:final} shows
that this value is related to
the $p$-rank of the class group of $\Q(\sqrt[p]{N})$.
As an illustration, when $N = 4261$, one computes that
the Hilbert class field
of $\Q(\sqrt[5]{4261})$ is $\Z/25\Z$. However, since the $5$-rank
of $\Z/25\Z$ is one, it follows that $e_5 = 1$.

\noindent \it Email addresses\rm:\tt \  fcale@math.harvard.edu

\hskip 22mm \tt \ emerton@math.northwestern.edu
\end{document}